\documentclass[final]{siamltex}

\usepackage{amsfonts,MnSymbol,graphicx}
\usepackage{multirow,booktabs}

\def\R{\mathbb{R}}
\def\RR{\mathbb{R}}
\def\wX{\widetilde X}

\def\RR{\mathbb{R}}

\newtheorem{Theorem}{Theorem}[section]
\newtheorem{Remark}[Theorem]{Remark}
\newtheorem{Lemma}[Theorem]{Lemma}
\newenvironment{Proof}{\begin{trivlist}\item[]
  {\em Proof.\ }}{$\square$\end{trivlist}}
\date{}
\begin{document}
\bibliographystyle{plain}

\title{Primal-dual interior-point multigrid method for\\
topology optimization\thanks{This work has been partly supported by Iraqi
Ministry of Higher Education and Scientific Research, Republic of Iraq and by
the EU FP7 project no.\ 313781 AMAZE.}}
\author{Michal Ko\v{c}vara\thanks{School of Mathematics, University of
    Birmingham, Edgbaston, Birmingham B15 2TT, UK, and Institute of Information Theory
and Automation, Academy of Sciences of the Czech Republic, Pod
vod\'arenskou v\v{e}\v{z}\'{\i}~4, 18208 Praha 8, Czech Republic}
\and Sudaba Mohammed\thanks{School of Mathematics, University of Birmingham, Edgbaston,
    Birmingham B15 2TT, UK, on leave from Department of Mathematics, University of Kirkuk, Iraq}}
\maketitle
%\tableofcontents
%\newpage

\begin{abstract}
An interior point method for the structural topology optimization is proposed.
The linear systems arising in the method are solved by the conjugate gradient
method preconditioned by geometric multigrid. The resulting method is then
compared with the so-called optimality condition method, an established
technique in topology optimization. This method is also equipped with the
multigrid preconditioned conjugate gradient algorithm. We conclude that, for
large scale problems, the interior point method with an inexact iterative
linear solver is superior to any other variant studied in the paper.
\end{abstract}



\section{Introduction}
The discipline of topology optimization offers challenging problems to
researchers working in large scale numerical optimization. The results are
essentially colors of pixels in a 2d or 3d ``pictures''. Hence, in order to
obtain high-quality results, i.e., fine pictures capturing all details, a very
large number of variables is essential. In this article we only consider the
discretized, finite dimensional topology optimization problem. For its
derivation and for general introduction to topology optimization, see, e.g.,
\cite{bendsoe-sigmund}.

We will consider the basic problem of topology optimization: minimization of
compliance under equilibrium equation constraints and the most basic linear
constraints on the design variables:
\begin{align}
  &\min_{x\in\R^m\!,\,u\in\R^n} \frac{1}{2}f^Tu\label{eq:to}\\
  &\mbox{subject to}\nonumber\\
  &\qquad K(x) u = f\nonumber\\
  &\qquad \sum_{i=1}^m x_i = V\nonumber\\
  &\qquad x_i\geq 0,\quad i=1,\ldots,m \nonumber\\
  &\qquad x_i\leq \overline{x},\quad i=1,\ldots,m\nonumber
\end{align}
where $K(x) = \sum_{i=1}^m x_i K_i$, $K_i\in\RR^{n\times n}$ and $f\in\RR^n$.
We assume that $K_i$ are symmetric and positive semidefinite and that
$\sum_{i=1}^m K_i$ is sparse and positive definite. We also assume that the
data $V\in\RR$ and $\overline{x}\in\RR^m$ is chosen such that the problem is
strictly feasible. For further reference, we will call the design variable $x$
the \emph{density}.

The most established and commonly used optimization methods to solve this
problem are the Optimality Conditions (OC) method
(\cite[p.308]{bendsoe-sigmund}) and the Method of Moving Asymptotes (MMA) by
Svanberg \cite{svanberg}. In both methods, the computational bottleneck
consists of the solution of a large scale linear system with a sparse symmetric
positive definite matrix (the equilibrium equation). This is traditionally used
by a direct solver, such as the Cholesky decomposition. Recently, several
authors proposed the use of iterative solvers, mostly preconditioned Krylov
subspace solvers, such as Conjugate Gradients (CG), MINRES or GMRES. These have
one big advantage which is specific for their use within optimization
algorithms: in the early (or even not-so-early) stages of the optimization
method, only a very low accuracy of the linear solver is needed. They also have
one big disadvantage: in the late stages of the optimization method, the linear
solvers become very ill-conditioned and thus a vanilla iterative method can
come into extreme difficulties.

It is therefore essential to use a good preconditioner for the Krylov subspace
method. The difficulty lies in the fact that as we approach the optimal
solution of the topology optimization problem, the condition number of the
stiffness matrices increases significantly. In fact, it is only controlled by
an artificial lower bound on the variable---if this bound was zero, the
stiffness matrix would be singular. Wang et al.\ \cite{wang} studied the
dependence of the condition number on the variables and concluded that it is a
combination of the ratio of maximum and minimum density and the conditioning of
a corresponding problem with constant density. Consequently, they proposed a
rescaling of the stiffness matrix combined with incomplete Cholesky
preconditioner. The rescaling results in constant order of condition number
during the optimization iterations. For large scale example still hundreds of
MINRES iterations are needed and hence the authors use recycling of certain
Krylov subspaces from previous iterations of the optimization method. Recently,
Amir et al.\ \cite{amir} proposed a multigrid preconditioner for the systems
resulting from OC or MMA methods and demonstrated that the resulting linear
system solver keeps its efficiency also for rapidly varying coefficient of the
underlying PDE, i.e., rapidly varying $x$ in (\ref{eq:to}).

While OC and MMA methods are the most popular methods in topology optimization,
they may not be the most efficient. The basic problem (\ref{eq:to}) is convex
(more precisely, it is equivalent to a convex problem) and we may thus expect
interior point methods to be highly efficient (see, e.g.,
\cite{nocedal-wright}). Indeed, Jarre et al.\ \cite{jarre-kocvara-zowe}
proposed an interior point method for the truss topology optimization problem
that is equivalent to the discretized problem (\ref{eq:to}), with the exception
that the stiffness matrix may be dense. They reported high efficiency of the
method and ability to solve large scale problems; they also proved convergence
of the proposed method. Maar and Schulz \cite{maar-schulz} studied interior
point methods for problem (\ref{eq:to}) with sparse stiffness matrices and
proposed to use a multigrid preconditioner for the GMRES method to solve the
arising indefinite linear systems.

A new comprehensive numerical study of optimization methods for topology
optimization can be found in \cite{rojas-stolpe}. The authors compare the
efficiency of different methods, including general purpose optimization solvers
such as SNOPT \cite{snopt}.

In this article we follow the path outlined by Jarre et al.\
\cite{jarre-kocvara-zowe} and by Maar and Schulz \cite{maar-schulz}. We use the
same interior point method as in \cite{jarre-kocvara-zowe} and, unlike in
\cite{maar-schulz}, reduce the linear systems to obtain positive definite
matrices. This allows us to use standard conjugate gradient method
preconditioned by standard V-cycle multigrid. We further use the same linear
solver in the OC method (in the same way as suggested in \cite{amir}) to get a
comparison with our interior point method. We will see that in both cases the
inexact multigrid preconditioned CG method leads to a very efficient
optimization solver. Most notably, in case of the interior point method we
obtain an approximately constant number of CG iterations needed to solve the
full problem which is independent of the size of the problem. In case of the OC
method, the total number of OC iterations is increasing with the problem size;
however, for a given problem size, the number of CG steps per one linear
systems remains almost constant, and very low, in all OC iterations,
notwithstanding the condition number of the stiffness matrix.

In this paper, we primarily consider the so-called variable thickness sheet
problem (\ref{eq:to}) and not its more popular cousin, the SIMP problem
\cite{bendsoe-sigmund}. The reason is the (hidden) convexity and existence of
solution of (\ref{eq:to}) (see \cite{ben1996hidden} and
\cite[p.272--274]{bendsoe-sigmund}). The goal of the paper is to study and
compare numerical methods for optimization problems. This can be done in a fair
way if the problem is convex; by introducing non-convexity, as in the SIMP
formulation, any such comparison is further influenced by many additional
factors. To demonstrate these difficulties and the fact that the iterative
solver is still a viable option in this context, we have added a brief section
on the SIMP model.

Finally, the methodology proposed in this paper is fully based on (typically
vectorizable and/or parallelizable) iterative schemes. It could thus be of
benefit to the \mbox{(re-)emerging} methods of distributed optimization
\cite{boyd2011distributed} and optimization on vector processors, namely GPU
\cite{schmidt20112589,zegard2013toward}, not only in the context of topology
optimization.

%----------------------------------------------------------------
\section{Newton systems for KKT conditions}
Let $\mu\in\R^n$, $\lambda\in\R$,
$\varphi\in\R^m$ and $\psi\in\R^m$ denote the respective Lagrangian multipliers
for constraints in (\ref{eq:to}).
The Karusch-Kuhn-Tucker (KKT) first order optimality conditions for (\ref{eq:to}) can be written as
\begin{align}
  -{\rm Res}^{(1)}:= &\ K(x) u - f = 0\label{eq:KKT1}\\
  -{\rm Res}^{(2)}:= &\ \sum_{i=1}^m x_i - V = 0\label{eq:KKT2}\\
  -{\rm Res}^{(3)}:= &\ -\frac{1}{2}u^TK_iu - \lambda - \varphi_i + \psi_i = 0,\quad i=1,\ldots,m\label{eq:KKT3}\\
  &\ \varphi_i x_i = 0,\quad i=1,\ldots,m\label{eq:KKT4}\\
  &\ \psi_i(\overline{x}-x_i)= 0,\quad i=1,\ldots,m\label{eq:KKT5}\\
  &\ x_i\geq 0,\quad \overline{x}-x_i\geq 0,\quad \varphi_i\geq 0,\quad \psi\geq 0\label{eq:KKT6}
\end{align}
We will perturb the complementarity constraints (\ref{eq:KKT4}) and
(\ref{eq:KKT5}) by barrier parameters $s,r>0$:
\begin{align}
  -{\rm Res}^{(4)}&:= \varphi_i x_i - s = 0,\quad i=1,\ldots,m\label{eq:KKT4p}\\
  -{\rm Res}^{(5)}&:= \psi_i(\overline{x}-x_i) - r= 0,\quad i=1,\ldots,m\label{eq:KKT5p}
\end{align}
and apply Newton's method to the system of nonlinear equations
(\ref{eq:KKT1}), (\ref{eq:KKT2}), (\ref{eq:KKT3}), (\ref{eq:KKT4p}),
(\ref{eq:KKT5p}).
In every step of the Newton method, we have to solve the linear system
\begin{equation}\label{eq:nwt}
  \begin{bmatrix} K(x) & 0 & B(u) & 0 & 0\\
    0 & 0 & e^T & 0 & 0\\
    B(u)^T & e & 0 & I & -I\\
    0 & 0 & \Phi & X & 0 \\
    0 & 0 & -\Psi & 0 & \wX \end{bmatrix}
    \begin{bmatrix} d_u\\d_\lambda\\d_x\\d_\varphi\\d_\psi\end{bmatrix}
    =
    \begin{bmatrix}{\rm Res}^{(1)}\\ {\rm Res}^{(2)}\\ {\rm Res}^{(3)}\\
    {\rm Res}^{(4)}\\ {\rm Res}^{(5)} \end{bmatrix}\,.
\end{equation}
Here $B(u) = (K_1u, K_2u,\ldots,K_mu)$, $e$ is a vector of all ones and
$$
  X={\rm diag}(x),\quad \wX={\rm diag}(\overline{x}-x),\quad
  \Phi = {\rm diag}(\varphi),\quad \Psi = {\rm diag}(\psi)
$$
are diagonal matrices with the corresponding vectors on the diagonal.

Because the last two equations only involve diagonal matrices, we can
eliminate $d_\varphi$ and~$d_\psi$:
\begin{align}
  d_\varphi &=X^{-1}({\rm Res}^{(4)} - \Phi d_t)\label{eq:phipsi1}\\
  d_\psi &=\wX^{-1}({\rm Res}^{(5)} - \Psi d_t)\,.\label{eq:phipsi2}
\end{align}
This will reduce the system (\ref{eq:nwt}) to
\begin{equation}\label{eq:nwtr}
  \begin{bmatrix} K(x) & 0 & B(u) \\
    0 & 0 & e^T \\
    B(u)^T & e & -(X^{-1}\Phi+\wX^{-1}\Psi) \end{bmatrix}
    \begin{bmatrix} d_u\\d_\lambda\\d_x\end{bmatrix}
    =
    \begin{bmatrix}{\rm Res}^{(1)}\\ {\rm Res}^{(2)}\\ \widetilde{\rm Res}^{(3)} \end{bmatrix}
\end{equation}
with
$$
  \widetilde{\rm Res}^{(3)} ={\rm Res}^{(3)}-X^{-1}{\rm Res}^{(4)}+\wX^{-1}{\rm Res}^{(5)}\,.
$$

We can now follow two strategies. Firstly, we can solve the system
(\ref{eq:nwtr}) as it is, i.e., an indefinite system of dimension
$m+n+1$. To simplify things, we can still eliminate the multipliers
$\varphi$ and $\psi$ as
$$
  \varphi_i = s/x_i,\quad \psi_i = r/(\overline{x}-x_i),\quad i=1,\ldots,m
$$
to get
\begin{equation}\label{eq:nwtr1}
  \begin{bmatrix} K(x) & 0 & B(u) \\
    0 & 0 & e^T \\
    B(u)^T & e & -(sX^{-2}+r\wX^{-2}) \end{bmatrix}
    \begin{bmatrix} d_u\\d_\lambda\\d_x\end{bmatrix}
    =
    \begin{bmatrix}{\rm Res}^{(1)}\\ {\rm Res}^{(2)}\\ \widetilde{\rm Res}^{(3)} \end{bmatrix}\,.
\end{equation}

\noindent{\bf Remark.} System (\ref{eq:nwtr1}) could be obtained
directly as a Newton system for optimality conditions of the following
``penalized'' problem:
\begin{align*}
  &\min_u \frac{1}{2}f^Tu + s\sum_{i=1}^m\log x_i + r\sum_{i=1}^m\log (\overline{x}-x_i)\\
  &\qquad{\rm s.t.}\quad K(x)u=f,\quad\sum_{i=1}^m x_i = V\,;
\end{align*}
see, e.g., \cite[Ch.19.1]{nocedal-wright}.

\medskip
Secondly, we can further reduce the Newton system (\ref{eq:nwtr}). As
the (3,3)-block matrix in (\ref{eq:nwtr}) is diagonal, we will compute
the Schur complement to the leading block to get
\begin{equation}\label{eq:nwtr2}
  Z \begin{bmatrix} d_u\\d_\lambda\end{bmatrix}  =
    {\rm Res}^{(Z)}    \,,
\end{equation}
with
\begin{equation}\label{eq:nwtr2Z}
  Z = \begin{bmatrix} K(x) & 0 \\
    0 & 0 \end{bmatrix} +  \begin{bmatrix}B(u) \\ e^T  \end{bmatrix}(X^{-1}\Phi+\wX^{-1}\Psi)^{-1}
    \begin{bmatrix}B(u)^T & e  \end{bmatrix}
\end{equation}
and
\begin{equation}\label{eq:nwtr2rhs}
    {\rm Res}^{(Z)} =
    \begin{bmatrix}{\rm Res}^{(1)}\\ {\rm Res}^{(2)} \end{bmatrix}
    + \begin{bmatrix}B(u) \\ e^T  \end{bmatrix}(X^{-1}\Phi+\wX^{-1}\Psi)^{-1}\widetilde{\rm Res}^{(3)}
    \,.
\end{equation}
The remaining part of the solution, $d_x$, is then computed by
\begin{equation}\label{eq:nwtr2Za}
  d_x = (X^{-1}\Phi+\wX^{-1}\Psi)^{-1}\left(\widetilde{\rm Res}^{(3)} -
  B^T{\rm Res}^{(1)} - e{\rm Res}^{(2)}\right)\,.
\end{equation}

%----------------------------------------------------------------
%----------------------------------------------------------------
%----------------------------------------------------------------
\section{Interior point method}
Once we have derived the Newton systems, the interior point algorithm is
straightforward (see, e.g., \cite[Ch.19]{nocedal-wright}). The details of the
individual steps of the algorithm will be given in subsequent paragraphs.

\subsection{The algorithm}
Denote $z=(u,\lambda,x,\varphi,\psi)^T$. Set $x_i=V/m,\ i=1,\ldots,m$, $u=K(x)^{-1}f$,  $\lambda=1,\ \varphi=e,\
\psi=e$. Set $s=1,\ r=1,\ \sigma_s, \sigma_r\in(0,1)$. Do until convergence:
\begin{enumerate}
\item Solve either system (\ref{eq:nwtr}) or (\ref{eq:nwtr2})  and
    compute the remaining components of vector $d$ from (\ref{eq:nwt}).
\item Find the step length $\alpha$.
\item Update the solution
$$
  z = z + \alpha d\,.
$$
\item If the stopping criterium for the Newton method is satisfied, update
    the barrier parameters
$$
  s = \sigma_s\cdot s,\quad r=\sigma_r\cdot r\,.
$$
Otherwise, keep current values of $s$ and $r$.\\
Return to Step~1.
\end{enumerate}

%--------------------------------------
\subsection{Barrier parameter update}
We use a fixed update of both parameters $s$ and $r$ with
$$
  \sigma_s = \sigma_r = 0.2 \,.
$$
This update leads to long steps and, consequently, small number of interior
point iterations. The value of the update parameter is a result of testing and
leads, in average, to the smallest overall number of Newton steps. A more
sophisticated version of the algorithm, with an adaptive choice of the barrier
parameters $s$ and $r$ can be found in \cite{jarre-kocvara-zowe}.
%----------------------------------------------------------------
\subsection{Step length}
We cannot take the full Newton step
$$
  z_{\rm new} = z + d
$$
because some variables could become infeasible with respect to the inequality
constraints (\ref{eq:KKT6}). We thus need to shorten the step in order to stay
strictly feasible with some ``buffer'' to the boundary of the feasible domain.
A simple step-length procedure is described below (see also
\cite[Ch.19.2]{nocedal-wright}).

Find $\alpha_l$ such that $x_i+(d_x)_i> 0$ for $i\in\{j:(d_x)_j<0\}$ and
$\alpha_u$ such that $x_i+(d_x)_i< \overline{x}$ for $i\in\{j:(d_x)_j>0\}$
using the following formulas:
$$
  \alpha_l = 0.9\cdot\min_{i:(d_x)_i<0}\left\{-\frac{x_i}{(d_x)_i}\right\},\quad
  \alpha_u = 0.9\cdot\min_{i:(d_x)_i>0}\left\{\frac{\overline{x}-x_i}{(d_x)_i}\right\}\,.
$$
The constant 0.9 guarantees the shortening of the step in the interior of the
feasible domain. Now take the smaller of these numbers and, if applicable,
reduce it to 1:
$$
  \alpha = \min\{\alpha_l,\alpha_u,1\}\,.
$$
A more sophisticated (and complicated) line-search procedure is
described in \cite{jarre-kocvara-zowe}.

It is worth noticing that for a properly chosen initial barrier parameter and
its update, the step-length reduction is almost never needed; this was, at
least, the case of our numerical examples and our choice of the parameters.
%----------------------------------------------------------------
\subsection{Stopping rules}
Following \cite{jarre-kocvara-zowe}, we terminate the Newton method whenever
$$
  \frac{\|{\rm Res}^{(1)}\|}{\|f\|} + \frac{\|\widetilde{\rm Res}^{(3)}\|}{\|\varphi\|+\|\psi\|} \leq \tau_{\scriptscriptstyle\rm NWT}.
$$

The full interior point method is stopped as soon as both parameters $s$ and
$r$ are smaller than a prescribed tolerance:
\begin{equation}\label{crit:ip}
  \max\{s,r\}\leq \tau_{\scriptscriptstyle\rm IP}\,.
\end{equation}

In our numerical experiments, we have used the values
$\tau_{\scriptscriptstyle\rm NWT}=10^{-1}$ and $\tau_{\scriptscriptstyle\rm
IP}=10^{-8}$.

\begin{Remark}
{\rm A more established criterium for terminating the interior point algorithm
would be to stop whenever all (scaled) residua are below some tolerance, i.e.,
$$
  \frac{\|{\rm Res}^{(1)}\|}{\|f\|} + \frac{\|\widetilde{\rm Res}^{(3)}\|}{\|\varphi\|+\|\psi\|}
   +\frac{\varphi^Tx}{\|\varphi\|\|x\|} +\frac{\psi^T(\overline{x}-x)}{\|\varphi\|\|x\|}
   \leq \tau_{\scriptscriptstyle\rm IP}.
$$
This criterium, however, leads to almost the same results as (\ref{crit:ip}),
hence we opted for the simpler and more predictable one.}
\end{Remark}

\begin{Remark}
{\rm The parameter $\tau_{\scriptscriptstyle\rm NWT}$ is kept constant in our
implementation, unlike in classic path-following methods. We will return to
this point later in Section~\ref{sec:exactIP}. }
\end{Remark}

%----------------------------------------------------------------
%----------------------------------------------------------------
%----------------------------------------------------------------
\section{Optimality Conditions method}
One of the goals of this paper is to compare the interior point method with the
established and commonly used Optimality Condition (OC) method. We will
therefore briefly introduce the basic algorithm and its new variant. For more
details, see (\cite[p.308]{bendsoe-sigmund}) and the references therein.

\subsection{OC algorithm}
Assume for the moment that the bound constraints in (\ref{eq:to}) are not
present. Then the KKT condition (\ref{eq:KKT3}) would read as
$$
  -u^T K_i u + \lambda = 0\,, \quad i=1,\ldots,m\,.
$$
(For convenience, we multiplied $\lambda$ from (\ref{eq:KKT3}) by $-\frac{1}{2}$.)
Multiplying both sides by $x_i$, we get
$$
  x_i\lambda = x_i u^T K_i u\,, \quad i=1,\ldots,m
$$
which leads to the following iterative scheme:
$$
  x_i^{\rm NEW} = \displaystyle\frac{1}{\lambda}\ x_i u^T K_i u\, , \quad i=1,\ldots,m\,.
$$
The new value of $x$ is then projected on the feasible set given by the bound
constraints. The value of $\lambda$ should be chosen such that $\sum_{i=1}^m
x_i^{\rm NEW} = V$ and is obtained by a simple bisection algorithm. Hence we
obtain the following algorithm called the OC method:
\paragraph{Algorithm OC}
Let $x\in\RR^m$ be given such that $\sum_{i=1}^m x_i = V$, $x\geq 0$. Repeat
until convergence:
\begin{enumerate}
\item $u=(K(x))^{-1} f$

\item $\overline{\lambda}=10000$, $\underline{\lambda}=0$

\item While $\overline{\lambda}-\underline{\lambda}>\tau_{\scriptscriptstyle \lambda}$

\begin{enumerate}
\item $\lambda = (\overline{\lambda}+\underline{\lambda})/2$
\item $x_i^{\rm NEW} = \min\left\{x_i \displaystyle\frac{u^T K_i u}{\lambda}\
    , \overline{x}\right\}\,, \quad i=1,\ldots,m$
\item $x = x^{\rm NEW}$
\item if $\sum_{i=1}^m x_i>V$ then set $\underline{\lambda}=\lambda$; else if $\sum_{i=1}^m x_i\leq V$ then set $\overline{\lambda}=\lambda$
\end{enumerate}
\end{enumerate}
The value of the bisection stopping criterium $\tau_{\scriptscriptstyle
\lambda}$ has been set to $10^{-11}$.

Notice that, due to positive semidefiniteness of $K_i$, the update in step 3(b)
is always non-negative and thus the lower-bound constraint in the original
problem (\ref{eq:to}) is automatically satisfied.

The basic version of the OC method converges (there are no known
counter examples) but is extremely slow. The reason for this is that,
from the very first iterations, the method is zig-zagging between two
clusters of points.
However, the following two modifications lead to a substantial improvement. To
the best of our knowledge, the second modification called Averaged OC is new.

\subsection{Damped OC}
\paragraph{Algorithm DOC}
Let $x\in\RR^m$ be given such that $\sum x_i = V$, $x\geq 0$. Repeat:
\begin{enumerate}
\item $u=(K(x))^{-1} f$

\item $x_i^{\rm NEW} = \min\left\{x_i \displaystyle\frac{(u^T K_i u)^q}{\lambda}\
    , \overline{x}\right\}\,, \quad i=1,\ldots,m$

\item $x = x^{\rm NEW}$
\end{enumerate}

Here $q$ is called the damping parameter; the typical choice is $q=1/2$. This
version of the method is widely used among the structural engineers.

\subsection{Averaged OC}
Let us define an operator $OC(\cdot)$ as a result of one step of the
standard OC algorithm.
\paragraph{Algorithm AOC}
Let $x\in\RR^m$ be given such that $\sum x_i = V$, $x\geq 0$. Repeat:
\begin{enumerate}
\item $x^{(1)} = OC(x)$

\item $x^{(2)} = OC(x^{(1)})$

\item $x = \frac{1}{2}(x^{(1)} + x^{(2)})$
\end{enumerate}
Numerical experiments suggest that Algorithm AOC is slightly faster than
Algorithm DOC. This modification seems to be new, at least we did not find it
in the existing literature.

\section{Multigrid conjugate gradient method}
In both optimization algorithms introduced above, we repeatedly need to solve
systems of linear equations. In this section, we will introduce an efficient
iterative method that seems to be most suitable for these problems. Throughout
this section, we assume that we want to solve the problem
\begin{equation}\label{eq:lineq}
  Az=b
\end{equation}
where $b\in\RR^n$ and $A$ is a $n\times n$ symmetric positive definite matrix.

\subsection{Multigrid method for linear systems}
Recall first the Correction Scheme (CS) version of the multigrid algorithm
(see, e.g., \cite{hackbusch}). Let $opt$ denote a (typically but not
necessarily) convergent iterative algorithm for (\ref{eq:lineq}):
$$
  z_{\rm new} = opt(A,b;z,\epsilon,\nu)\,,
$$
where, on input, $z$ is the initial approximation of the solution, $\epsilon$
is the required precision and $\nu$ the maximum number of iterations allowed.
This will be called the \emph{smoother}. A typical example is the Gauss-Seidel
iterative method.

Assume that there exist $\ell$ linear operators
$I_k^{k-1}:\RR^{n_k}\to\RR^{n_{k-1}}$, $k=2,\ldots,\ell$, with
$n:=n_\ell>n_{\ell-1}>\cdots>n_2>n_1$ and let $I_{k-1}^k:=(I_k^{k-1})^T$. These
are either constructed from finite element or finite difference refinements of
some original coarse grid (geometric multigrid) or from the matrix $A$
(algebraic multigrid); see \cite{briggs2000multigrid} for details.

Define the ``coarse level'' problems
$$
  A_k z_k=b_k,\quad k=1,\ldots,\ell-1
$$
with
$$
  A_{k-1} = I_k^{k-1} (A_k) I_{k-1}^k,\quad b_{k-1} = I_k^{k-1} (b_k),\quad k=2,\ldots,\ell\,.
$$

%----------------------------------------------------------------
\paragraph{Algorithm MG} (V-cycle correction scheme multigrid)

\begin{description}
\item Set $\epsilon,\epsilon_0$. Initialize $z^{(\ell)}$.
\item for $i=1:niter$
\begin{description}
\item $z^{(\ell)} := mgm(\ell,z^{(\ell)},b_\ell)$
\item test convergence
\end{description}
\item end

\item function
    $z^{(k)}=mgm(k,z^{(k)},r_k)$
\begin{description}
\item if $k=0$
\begin{description}\item $z^{(k)}:= opt(A_1,b_1;z^{(k)},\epsilon_0,\nu_0)$\hfill (coarsest grid solution)  \end{description} %
\item else%
\begin{description}\item $z^{(k)}:= opt(A_k,b_k;z^{(k)},\epsilon,\nu_1)$ \hfill (pre-smoothing)
\item $r_{k-1} = I_k^{k-1} (r_k - A_k z^{(k)})$\hfill (restricted residuum)
\item $v^{(k-1)} =
    mgm(k-1,0_{n_{k-1}},r_{k-1})$\hfill (coarse grid correction)
\item $z^{(k)} := z^{(k)} + I_{k-1}^k v^{(k-1)}$\hfill (solution update)
\item $z^{(k)}:= opt(A_k,b_k;z^{(k)},\epsilon,\nu_2)$  \hfill (post-smoothing) \end{description}%
\item end
\end{description}
\end{description}

\subsection{Multigrid preconditioned conjugate gradient method}
Although the multigrid method described above is very efficient, an even more
efficient tool for solving  (\ref{eq:lineq}) may be the preconditioned
conjugate gradient (CG) method, whereas the preconditioner consist of one step
of the V-cycle multigrid method. The algorithm is described below (see, e.g.,
\cite{golub-vanloan}).

\paragraph{Algorithm PCG}
\begin{description}
\item Given initial $z$, set $r := Az-b$
\item $y := mgm(\ell,0_{n},r)$
\item Set $p:=-y$

\item for $i=1:niter$
\begin{description}
\item $\alpha:=\displaystyle\frac{r^T y}{p^TAp}$
\item $z:= z + \alpha p$
\item $\tilde{r} := r + \alpha Ap$
\item $\tilde{y} := mgm(\ell,0_{n},\tilde{r})$
\item $\beta := \displaystyle\frac{\tilde{r}^T\,\tilde{y}}{r^T\, y}$
\item $p:= -y +\beta p$
\item $r:= \tilde{r}$, $y:= \tilde{y}$
\item test convergence
\end{description}
\item end
\end{description}

%----------------------------------------------------------------
\section{Multigrid conjugate gradients for IP and OC methods}
The main goal of this section (and of the whole article) is to study the effect
of the multigrid preconditioned CG method in the IP and OC algorithms. We will
also compare them to their counterparts, IP and OC with direct solvers.

The details on discretization and the choice of prolongation and restriction
operators will be given in Section~\ref{sec:examples}.

\subsection{Multigrid conjugate gradients for IP}
Our goal is to solve the linear systems arising in the Newton method, by the
conjugate gradient method preconditioned by one V-type multigrid step. We can
choose one of the three equivalent systems to solve, namely the full system
(\ref{eq:nwt}), the reduced saddle-point system (\ref{eq:nwtr}) and the
so-called augmented system (\ref{eq:nwtr2}). We prefer the last one for the
following reasons.
\begin{itemize}
\item The matrix $Z$ in (\ref{eq:nwtr2}) is positive definite and we can thus readily apply the standard conjugate gradient method together with the standard V-cycle as a preconditioner. We could, of  course, use GMRES or MINRES for the indefinite systems in (\ref{eq:nwt}) and (\ref{eq:nwtr}), however, the multigrid preconditioner, in particular the smoother, would become more complicated in this case; see \cite{maar-schulz}, who used so-called transforming smoothers introduced by Wittum \cite{wittum}.
\item In order to use the multigrid preconditioner, we have to define prolongation/restriction operators for the involved variables. This can be easily done in case of the system (\ref{eq:nwtr2}) that only involves the displacement variable $u\in\RR^n$ plus one additional variable $\lambda$, the Lagrangian multiplier associated with the volume constraint; see the next Section~\ref{sec:examples} for details.

    If, on the other hand, we decided to solve (\ref{eq:nwt}) or (\ref{eq:nwtr}), we would have to select an additional restriction operator for the variables associated with the finite elements; this operator should then be ``compatible'' with the nodal-based restriction operator. This is a rather non-trivial task and can be simply avoided by choosing system (\ref{eq:nwtr2}).
\end{itemize}

\begin{figure}[h]
\begin{center}
 \resizebox{0.45\hsize}{!}
   {\includegraphics{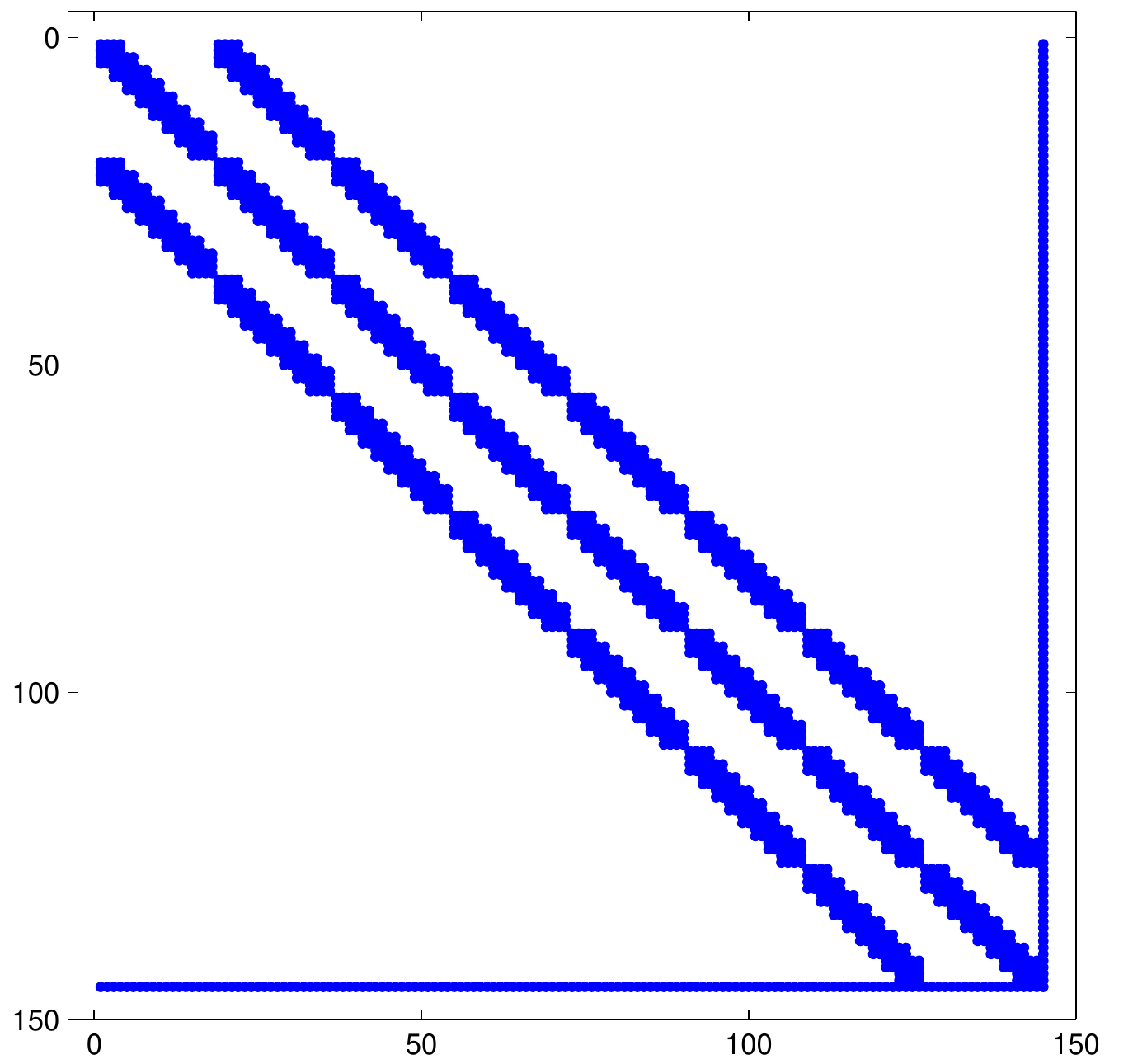}}\quad
 \resizebox{0.45\hsize}{!}
   {\includegraphics{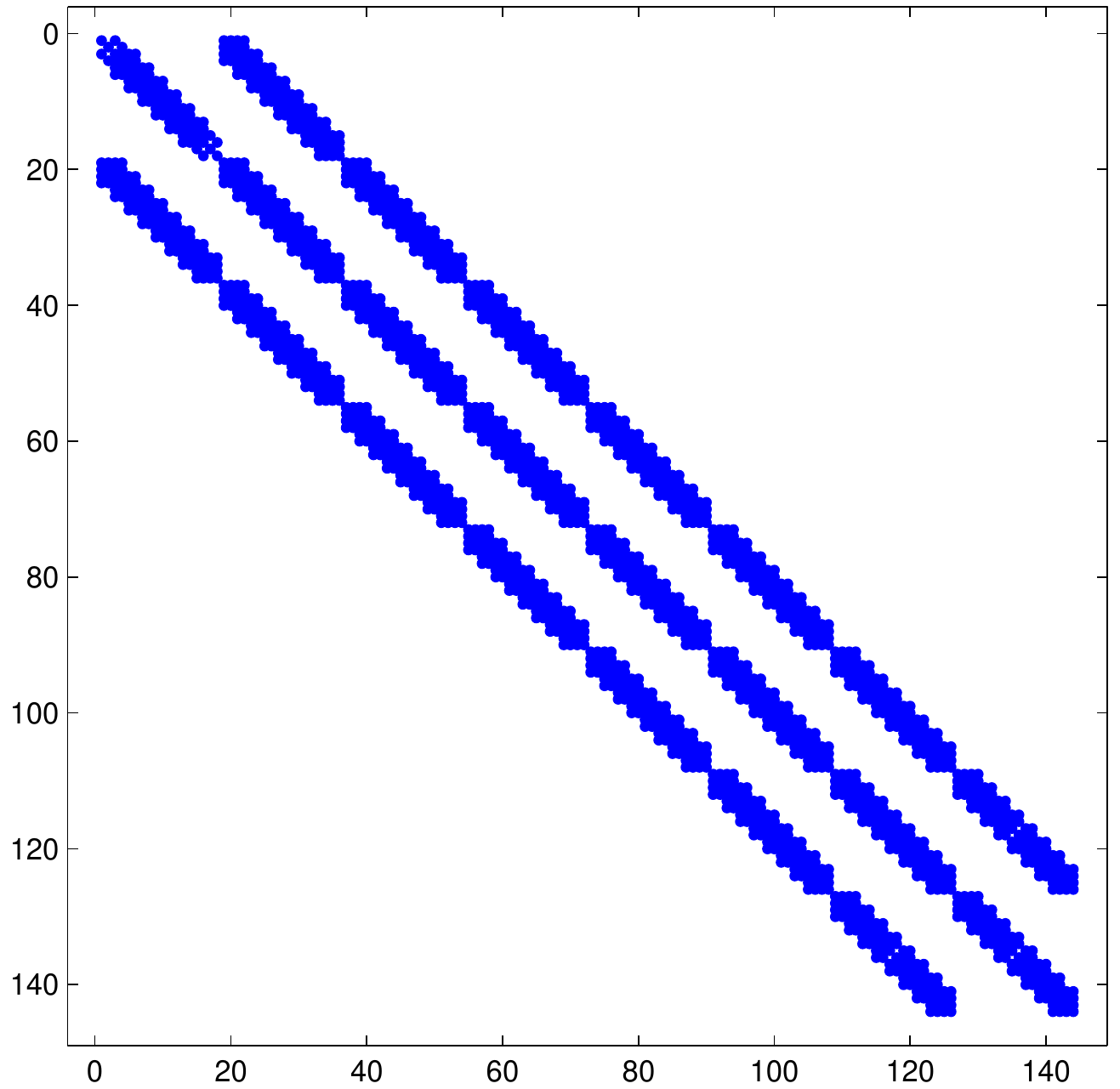}}
  \end{center}
  \caption{Typical sparsity structure of matrix $Z$ from the augmented system (\ref{eq:nwtr2}) (left) and of the stiffness matrix $K$ (right)}\label{fig:matrix}
\end{figure}

The matrix $Z$ from (\ref{eq:nwtr2}) is positive definite, sparse and typically
has an arrow-type sparsity structure: it is banded apart from the last full row
and column; see Figure~\ref{fig:matrix}-left. The bandwidth grows,
approximately, with the square root of the problem size. At the same time, the
number of non-zeros in each row is always the same, notwithstanding the problem
size.
%%%%%%%%%%%%%%%%%%%%%%%%%%%%%%%%%%%%%%%%
%%%%%%%%%%%%%%%%%%%%%%%%%%%%%%%%%%%%%%%%

%----------------------------------------------------------------
\paragraph{Stopping rule}
It is a big advantage of iterative methods, over direct solvers, that they
allow us to control the precision of the approximate solution and stop whenever
even a low required precision is reached. In our implementation, the PCG method
is stopped whenever
\begin{equation}\label{eq:cgstop}
  \|\rho\|\,\|b\|\leq 10^{-2}
\end{equation}
where $\rho$ is residuum and $b$ the right-hand side of the linear system,
respectively. In this way we only compute an approximate Newton direction; it
is shown, e.g., in \cite{dembo} that the resulting method converges once the
approximate Newton direction is ``close enough'' (though not infinitesimally
close in the limit) to the exact solution of the Newton system. Furthermore,
for convex quadratic programming problems, Gondzio \cite{gondzio} has shown
that when the PCG method is stopped as soon as $\|\rho\|\leq 0.05 s$ ($s$ being
the barrier parameter), the theoretical complexity of the interior point method
is the same as with the exact linear solver. Inexact iterative solvers in the
context of other optimization problems and algorithms were further studied,
e.g., in \cite{conn-gould-toint,pennon-iter,mizuno-jarre,toh}.

In our case, the value of $10^{-2}$ proved to be a good compromise between the
overall number of Newton steps and the overall number of PCG iterations within
the IP method. With this stopping criterium, the IP methods requires,
typically, 2--4 PCG iterations in the initial and in many subsequent IP steps.
Only when we get close to the required accuracy, in the last 2--3 IP steps, the
conditioning of the matrix $Z$ increases significantly and so does the number
of PCG steps, typically to 10--30; see the next section for detailed numerical
results.

\subsection{Multigrid conjugate gradients for OC}\label{sec:CGOC}
Within the OC algorithm, the multigrid CG method will be used to solve the
discretized equilibrium equation $Ku=f$. Recall that $K$ is assumed to be a
positive definite matrix. Moreover $K$ is very sparse and, if a reasonably good
numbering of the nodes is used, banded. A typical non-zero structure of $K$ is
shown in Figure~\ref{fig:matrix}-right: it is exactly the same as for the
matrix in (\ref{eq:nwtr2}) in the IP method, apart from the additional last
column and row in the augmented matrix in (\ref{eq:nwtr2}).
%We will return to this subtle difference in the next section.

%\subsection{OC stopping criteria}
The only degrees of freedom in the resulting algorithm are the stopping
criteria for the OC method and for the multigrid CG method.

\paragraph{The overall stopping criterium}
As the dual information (Lagrangian multipliers associated with the bound
constraints) is not readily available, so far the only practical (and widely
used) stopping criterium for the OC method is the difference in the objective
function value in two subsequent iterations. Needless to say that, unless we
have an estimate for the rate of convergence, this criterium can be misleading
and may terminate the iteration process long before some expected approximation
of the optimum has been reached. Nevertheless, many numerical experiments
suggest that this criterium is not as bad as it seems and serves its purpose
for the OC method.

Hence the OC method is typically stopped as soon as
\begin{equation}\label{eq:OCstop}
  |f^T u_k - f^T u_{k-1}| \leq \tau_{\scriptscriptstyle\rm OC}
\end{equation}
where $k$ is the iteration index. In our numerical experiments we have used
$\tau_{\scriptscriptstyle\rm OC}=10^{-5}$; this value has been chosen such that
the OC results are comparable to the IP results, in the number of valid digits
both in the objective function and in the variables; see
Section~\ref{sec:exact} for more details.

\paragraph{Stopping criterium for the multigrid CG method}
As already mentioned above, one of the advantages of an iterative method is the
fact that an exact solution to the linear system is not always needed. In such
a case, we can stop the iterative method after reaching a relatively low
accuracy solution. The required accuracy of these solutions (such that the
overall convergence is maintained) is well documented and theoretically
supported in case of the IP method; it is, however an unknown in case of the OC
method; see \cite{amir} for detailed discussion. Clearly, if the linear systems
in the OC method are solved too inaccurately, the whole method may diverge or
just oscillate around a point away from the solution.

We have opted for the following heuristics that guarantees the (assumed)
overall convergence of the OC method. Notice that the OC method is a feasible
descent algorithm. That means that every iteration is feasible and the
objective function value in the $k$-th iteration is smaller than that in the
$(k-1)$-st iteration. Hence
\begin{itemize}
\item we start with $\tau=10^{-4}$\,;
\item if $f^T u_k > f^T u_{k-1}$, we update $\tau := 0.1\, \tau$.
\end{itemize}
In our numerical tests, the update had to be done only in few cases and the
smallest value of $\tau$ needed was  $\tau=10^{-6}$. Recall that this is due to
our relatively mild overall stopping criterium (\ref{eq:OCstop}). In the next
section, we will see that this heuristics serves its purpose, as the number of
OC iterations is almost always the same, whether we use an iterative or a
direct solver for the linear systems.

\section{Numerical experiments}\label{sec:examples}
This section contains detailed results of three numerical examples. All codes
were written entirely in MATLAB. Notice, however, that when we refer to a
direct solver for the solution of linear system, we mean the backslash
operation in MATLAB which, for our symmetric positive definite systems, calls
the CHOLMOD implementation of the Cholesky method \cite{cholmod}. This
implementation is highly tuned, very efficient and written in the C language.
So whenever we compare CPU times of the iterative solver with the direct
solver, we should keep this in mind. These comparisons are given solely to show
the tendency in the CPU time when increasing the problem size. All problems
were solved on an Intel Core i5-3570 CPU at 3.4GHz with 8GB RAM, using MATLAB
version 8.0.0 (2012b) running in 64 bit Windows 7.

In all examples, we use square finite elements with bilinear basis functions
for the displacement variable $u$ and constant basis functions for the
thickness variable $x$, as it is standard in topology optimization. The
prolongation operators $I_{k-1}^k$ for the variable $u$ are based on the
nine-point interpolation scheme defined by the stencil $\begin{pmatrix}
\frac{1}{4}&\frac{1}{2}&\frac{1}{4}\\
\frac{1}{2}&1&\frac{1}{2}\\
\frac{1}{4}&\frac{1}{2}&\frac{1}{4}\\
\end{pmatrix}$; see, e.g., \cite{hackbusch}.
When solving the linear system (\ref{eq:nwtr2}) in the interior point method,
we also need to prolong and restrict the single additional variable $\lambda$;
here we simply use the identity.

The examples are solved with isotropic material with Young's modulus equal to 1
and Poisson's ratio 0.3. The physical dimensions of the computational domain
are given by the coarsest mesh, whereas the coarse level element has dimension
$1\times 1$. The upper bound on the variable $x$ is set to $\overline{x}=2$.
The load is always defined on the finest discretization level on edges of two
elements sharing a node on the boundary specified in each example. The load
always acts in vertical direction. Thus the non-zero elements of the
discretized load vector will be $(-\frac{1}{2}, -1, -\frac{1}{2})$, associated
with the vertical components of the specified boundary nodes its two immediate
neighbours on this boundary.

The meaning of the captions in the following tables:
\begin{description}
\setlength\itemsep{0em}
\item[] {\bf problem}\ldots{}the first two numbers describe the dimension
    of the computational domain, the last number is the number of mesh
    refinements
\item[] {\bf variables}\ldots{}number of variables in the linear systems
\item[] {\bf feval}\ldots{}total number of function evaluations (equal to
    the number of linear systems solved)
\item[] {\bf total CG iters}\ldots{}total number of CG iterations in the optimization process
\item[] {\bf solver CPU time}\ldots{}total CPU time spent in the solution of linear systems
\item[] {\bf average CG iters}\ldots{}average number of CG iterations per one linear system
%\item[] {\bf complexity factor}\ldots{}assuming that the computational complexity of the linear system solver is polynomial in the problems dimension $n$, we compute the approximate degree $d$ of the polynomial $cn^d$ from the results of two subsequent problem refinements; so, e.g., for complexity factor equal to one, we get linear complexity, for 2 quadratic complexity, etc. If the factor is smaller than one, we obtain sub-linear complexity.
\end{description}

%\paragraph{Remark} Notice that the complexity factor introduced above may not correspond to the theoretical complexity of the method, as given, e.g., in \cite{golub-vanloan}. Our matrices are banded and the bandwidth grows, approximately, with the square root of the problem size. At the same time, the number of non-zeros in each row is always the same, notwithstanding the problem size. Moreover, the iterative method only solves the problems approximately and is stopped long before finding a solution with the same precision as the direct method.

\subsection{Example 1}
\begin{figure}[h]
\begin{center}
 \resizebox{0.4\hsize}{!}
   {\includegraphics{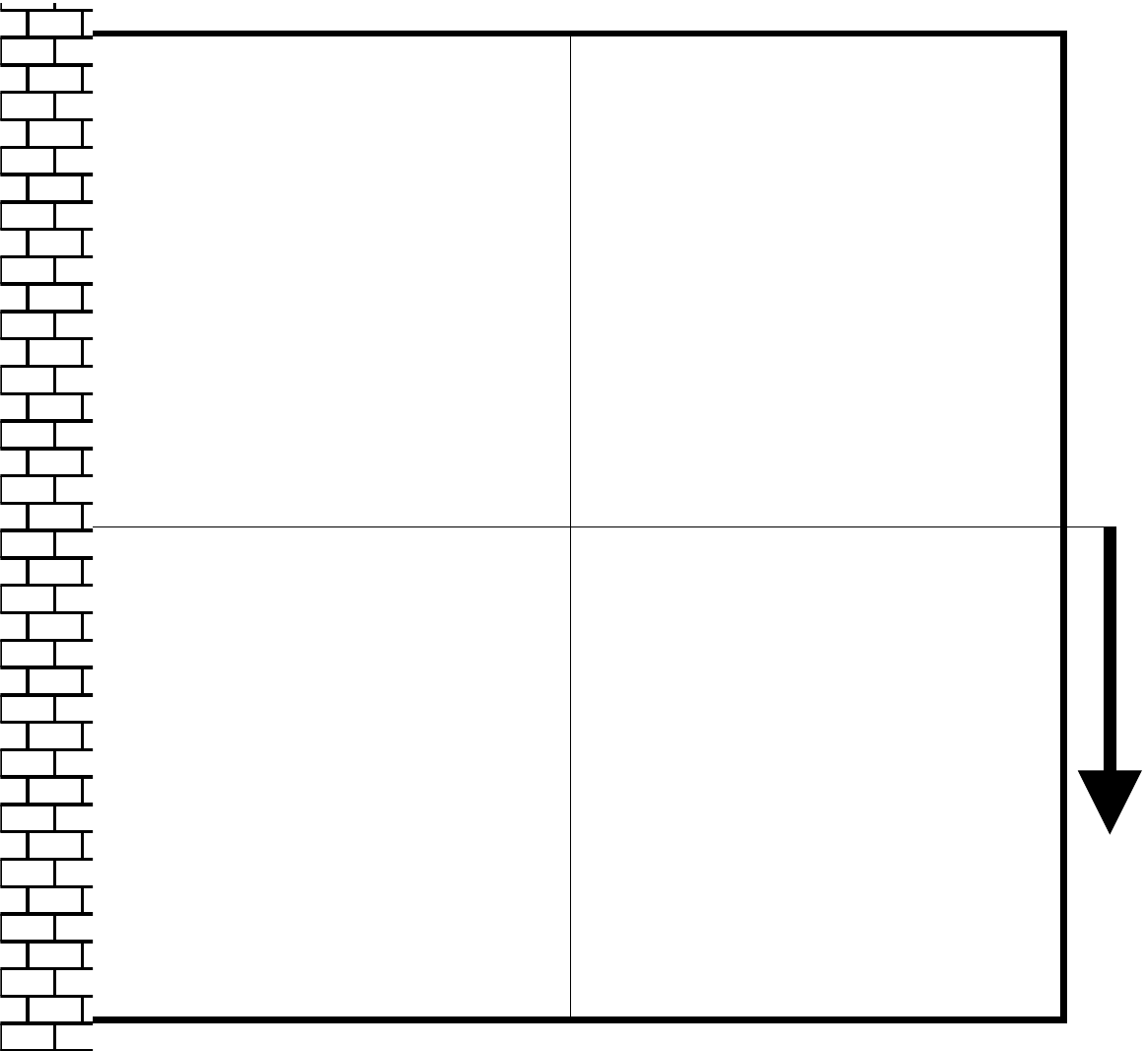}}\quad
    \resizebox{0.4\hsize}{!}
   {\includegraphics{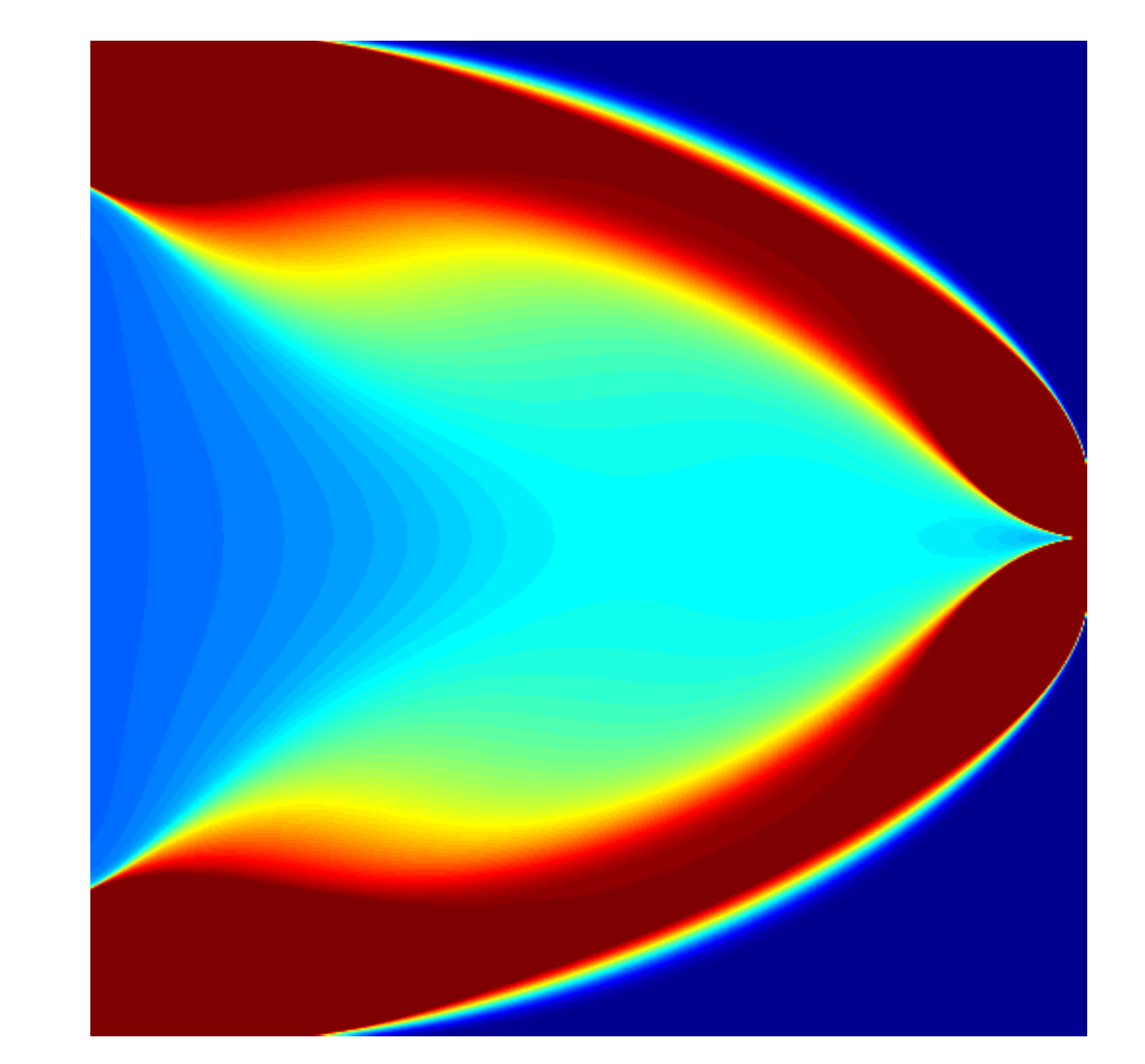}}
  \end{center}
  \caption{Example~1, initial setting with coarsest mesh and optimal solution.}\label{fig:11}
\end{figure}

We consider a square computational domain with the coarsest mesh consisting of
$2\times 2$ elements. All nodes on the left-hand side are fixed, the right-hand
middle node is subject to a vertical load; see Figure~\ref{fig:11}. We use up
to nine refinements levels with the finest mesh having 262\,144 elements and
525\,312 nodal variables (after elimination of the fixed nodes).
%%%%%%%%%%%%%%%%%%%%%%%%%%%%%%%%%%%%%%%%
%%%%%%%%%%%%%%%%%%%%%%%%%%%%%%%%%%%%%%%%

Table~\ref{tab:1} presents the results of the interior point method. We can see
that, with increasing size of the problem, the total number of CG iterations is
actually decreasing. This is due to our specific stopping criterium explained
in the previous section. We also observe that the average number of CG
iterations per linear system is very low and, in particular, \emph{is not
increasing with the problem size}, the result of the multigrid preconditioner.
% Table generated by Excel2LaTeX from sheet '22_school (2)'
\begin{table}[htbp]
  \centering
  \caption{Example 1, interior point method with iterative solver}
    \begin{tabular}{crcccc}
    \toprule
              &       &       &  total     &   solver    & average \\
    problem & variables  & feval & CG iters & CPU time & CG iters \\
    \midrule
    223   & 145      & 31    & 253   & 0.18  & 8.16 \\
    224   & 545      & 30    & 281   & 0.44  & 9.37 \\
    225   & 2\,113   & 29    & 197   & 0.91  & 6.79 \\
    226   & 8\,321   & 28    & 139   & 2.79  & 4.96 \\
    227   & 33\,025  & 27    & 119   & 12.7  & 4.41 \\
    228   & 131\,585 & 25    & 104   & 45.8  & 4.16 \\
    229   & 525\,313 & 27    &  85   & 156.0 & 3.15 \\
    \bottomrule
    \end{tabular}%
  \label{tab:1}%
\end{table}%

Let us now compare these results with those for the OC method where the linear
system is just the equilibrium equation; see Table~\ref{tab:2}. As expected,
the number of OC iterations (and thus the number of linear systems and the
total number of CG iterations) \emph{grows} with the size of the problem. Also
in this case the average number of CG iterations is almost constant,
notwithstanding the size of the problem.
% Table generated by Excel2LaTeX from sheet '22_school (2)'
\begin{table}[htbp]
  \centering
  \caption{Example 1, OC method with iterative solver}
    \begin{tabular}{crcccc}
    \toprule
              &       &       &  total     &   solver    & average  \\
    problem & variables  & feval & CG iters & CPU time & CG iters  \\
    \midrule
    223   & 144       & 19    & 56    & 0.04  & 2.95 \\
    224   & 544       & 33    & 100   & 0.14  & 3.03 \\
    225   & 2\,112    & 55    & 164   & 0.65  & 2.98 \\
    226   & 8\,320    & 85    & 254   & 4.84  & 2.99 \\
    227   & 33\,024   & 111   & 332   & 30.8  & 2.99 \\
    228   & 131\,584  & 119   & 362   & 133.0 & 3.04 \\
    229   & 525\,312  & 123   & 368   & 636.0 & 2.99 \\
    \bottomrule
    \end{tabular}%
  \label{tab:2}%
\end{table}%

\begin{figure}[h]
\begin{center}
 \resizebox{0.50\hsize}{!}
   {\includegraphics{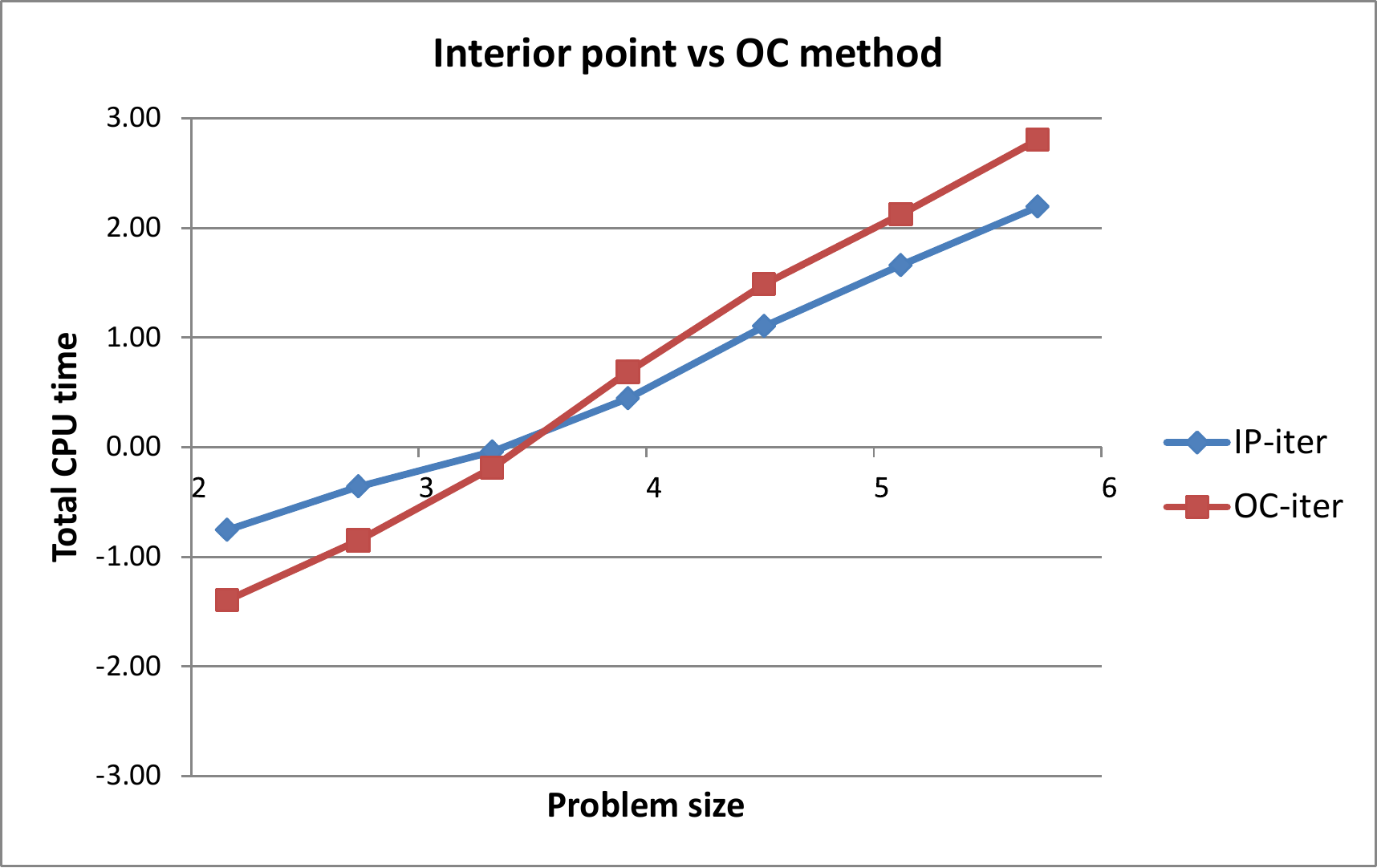}}\
 \resizebox{0.46\hsize}{!}
   {\includegraphics{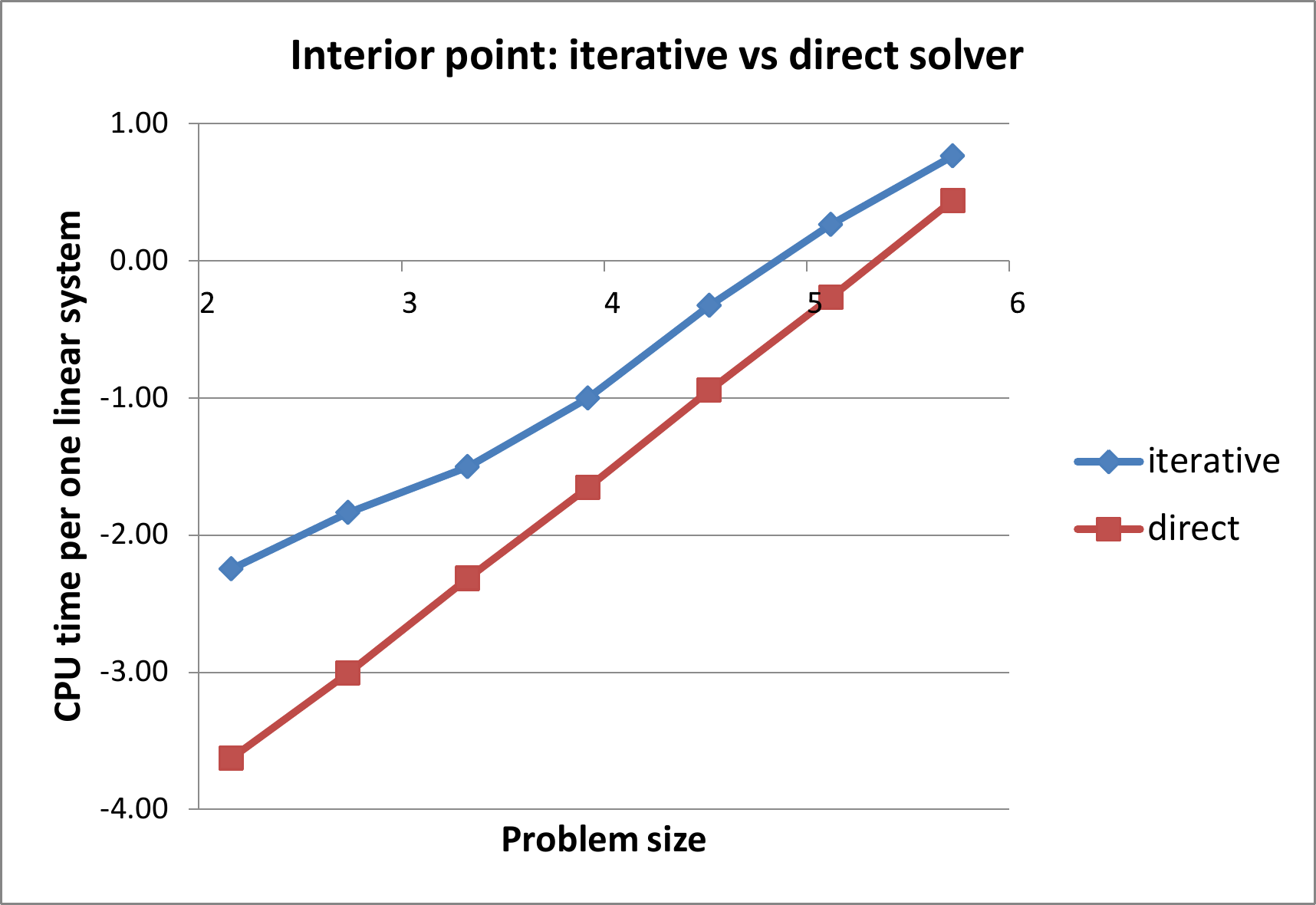}}
  \end{center}
  \caption{Example~1, \emph{left}: total CPU time spent in the iterative linear solver for the interior point and the OC method; \emph{right}: interior point method, total CPU time spent in the iterative linear solver and in the direct solver.}\label{fig:12}
\end{figure}

The comparison of the interior point method with the OC method is graphically
presented in Figure~\ref{fig:12} (left). Here we can see, in the log-log scale,
the total CPU time spent in the linear solver, growing with the size of the
problem. While initially worse than the OC method, the interior point method
grows slower and soon catches up and overtakes the OC method. For both methods,
the growth is almost linear for the larger problems, so that we can estimate
the growth in the CPU time as a polynomial function $cn^d$ of the problem
dimension $n$. For the interior point method, the degree $d=0.907$ while for
the OC method $d=1.09$. This means that the overall computational complexity of
the IP method with inexact Newton and inexact multigrid CG methods is slightly
\emph{sublinear}. For the OC method, it is just a bit worse than linear.
%%%%%%%%%%%%%%%%%%%%%%%%%%%%%%%%%%%%%%%%
%%%%%%%%%%%%%%%%%%%%%%%%%%%%%%%%%%%%%%%%

In Figure~\ref{fig:12} (right) we compare the iterative solver used in the
interior point method with a direct Cholesky solver (see the warning at
beginning of this section!). We can clearly see that the time for the (C coded)
direct solver grows quicker than for the (MATLAB coded) iterative solver.

\subsection{Example 2}
\begin{figure}[h]
\begin{center}
 \resizebox{0.48\hsize}{!}
   {\includegraphics{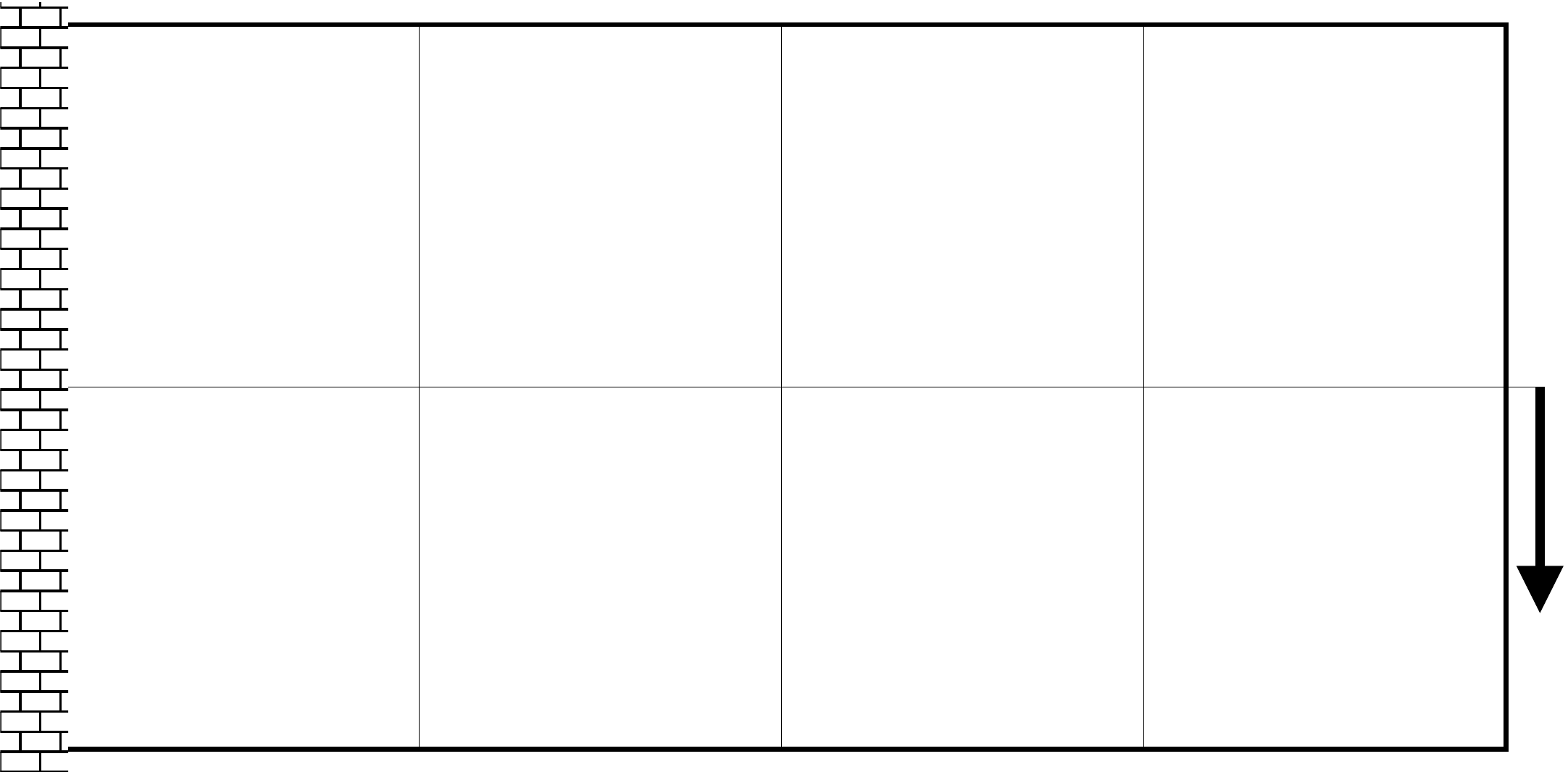}}\
    \resizebox{0.48\hsize}{!}
   {\includegraphics{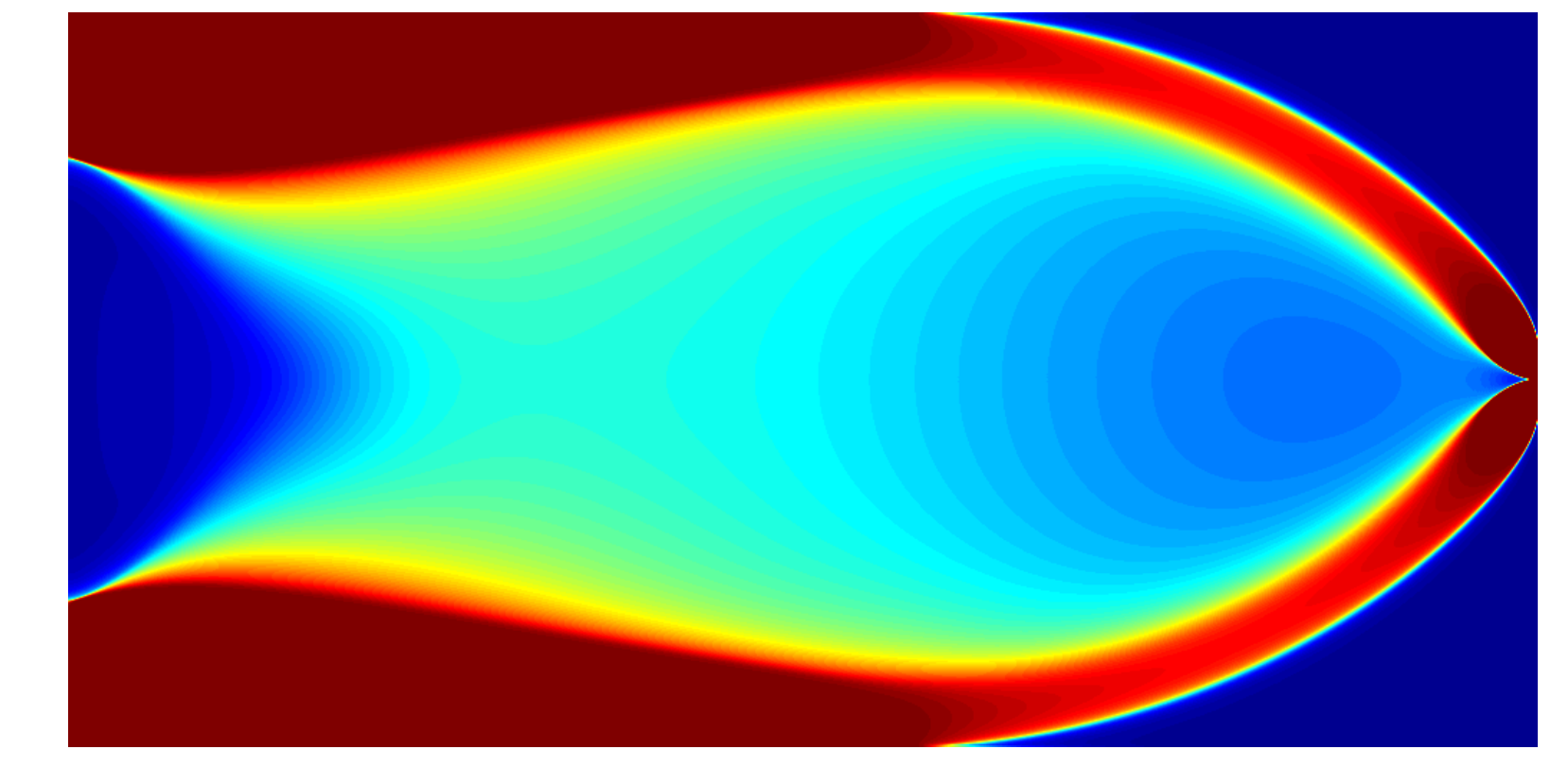}}
  \end{center}
  \caption{Example~2, initial setting with coarsest mesh and optimal solution.}\label{fig:21}
\end{figure}

The next example is similar to the previous one, only the computational domain
is ``longer'' in the horizontal direction; the coarsest mesh consists of
$4\times 2$ elements. It is well known that the conditioning of this kind of
examples grows with the slenderness of the domain. As before, all nodes on the
left-hand side are fixed, the right-hand middle node is subject to a vertical
load; see Figure~\ref{fig:21}. Again, we use up to nine refinements levels with
the finest mesh having 524\,288 elements and 1\,050\,624 nodal variables (after
elimination of the fixed nodes).
%%%%%%%%%%%%%%%%%%%%%%%%%%%%%%%%%%%%%%%%
%%%%%%%%%%%%%%%%%%%%%%%%%%%%%%%%%%%%%%%%

We first show the results of the interior point method in Table~\ref{tab:3}.
Just as in the previous example, the total number of CG iterations is
decreasing with the increasing size of the problem. Again, the average number
of CG iterations per linear system is very low and not increasing.
% Table generated by Excel2LaTeX from sheet '42_school (2)'
\begin{table}[htbp]
  \centering
  \caption{Example 2, IP method with iterative solver}
    \begin{tabular}{crcccc}
    \toprule
              &       &       &  total     &   solver    & average  \\
    problem & variables  & feval & CG iters & CPU time & CG iters  \\
    \midrule
    423   & 288         & 33    & 265   & 0.24  & 8.03 \\
    424   & 1\,088      & 32    & 342   & 0.87  & 10.69\\
    425   & 4\,224      & 31    & 207   & 1.89  & 6.68 \\
    426   & 16\,640     & 30    & 160   & 7.77  & 5.33 \\
    427   & 66\,048     & 29    & 139   & 31.1  & 4.79 \\
    428   & 263\,168    & 27    & 123   & 119.0 & 4.56 \\
    429   & 1\,050\,624 & 27    & 101   & 385.0 & 3.74 \\
    \bottomrule
    \end{tabular}%
  \label{tab:3}%
\end{table}%
Compare this with the OC solver results in Table~\ref{tab:4}. In this case, we
only consider eight refinement levels, as the largest problem would take too
much time on our computer. Contrary to the previous example, the average number
of CG iterations is slightly increasing due to the worse conditioning.
% Table generated by Excel2LaTeX from sheet '42_school (2)'
\begin{table}[htbp]
  \centering
  \caption{Example 2, OC method with iterative solver}
    \begin{tabular}{crccccc}
    \toprule
              &       &       &  total     &   solver    & average  \\
    problem & variables  & feval & CG iters & CPU time & CG iters  \\
    \midrule
    423   & 288      & 39   & 117   & 0.13   & 3.00 \\
    424   & 1\,088   & 45   & 144   & 0.34   & 3.20 \\
    425   & 4\,224   & 77   & 262   & 2.10   & 3.40 \\
    426   & 16\,640  & 123  & 423   & 16.2   & 3.44 \\
    427   & 66\,048  & 157  & 542   & 97.1   & 3.45 \\
    428   & 263\,168 & 165  & 739   & 552    & 4.48 \\
    \bottomrule
    \end{tabular}%
  \label{tab:4}%
\end{table}%

\begin{figure}[h]
\begin{center}
 \resizebox{0.50\hsize}{!}
   {\includegraphics{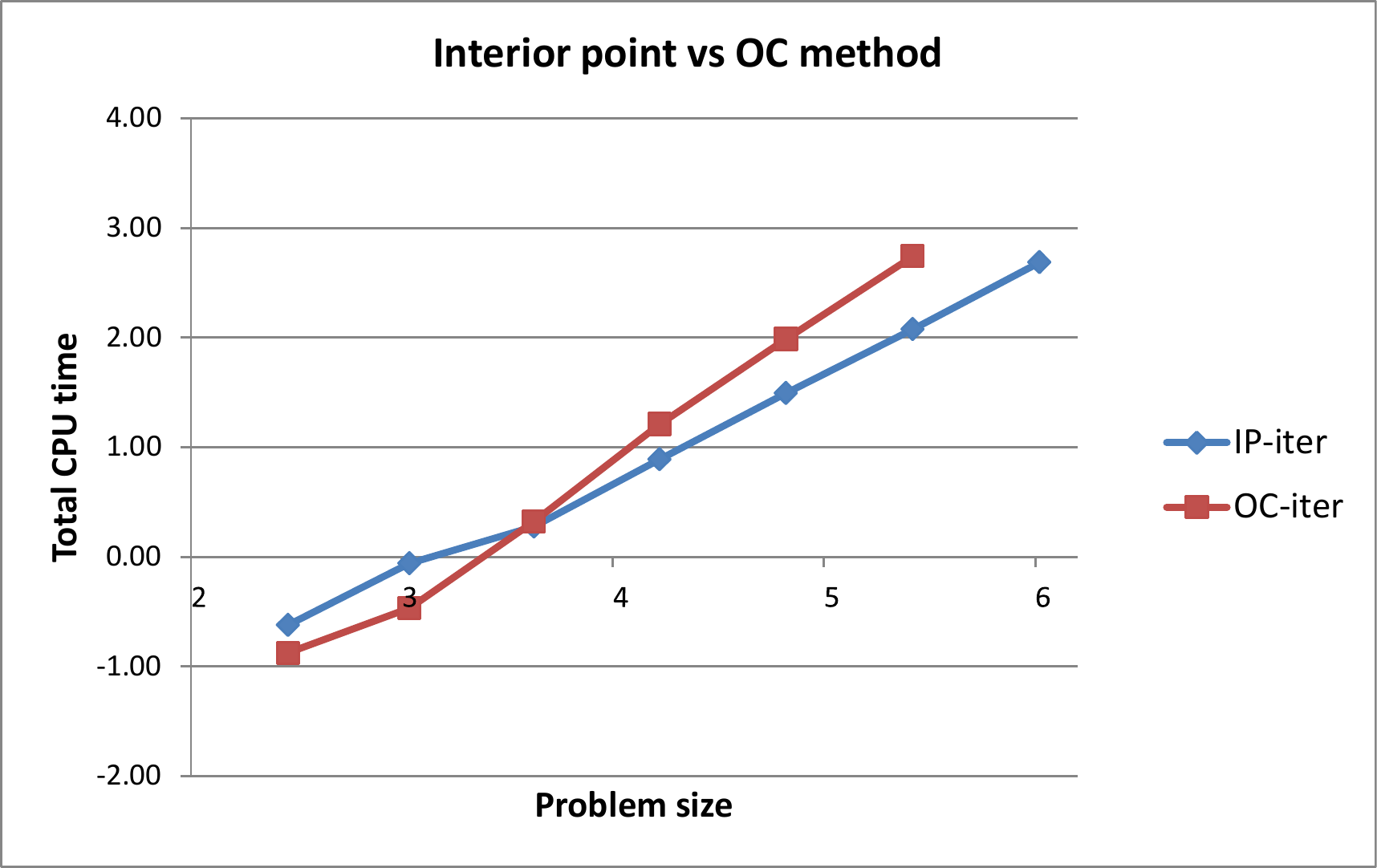}}\
 \resizebox{0.46\hsize}{!}
   {\includegraphics{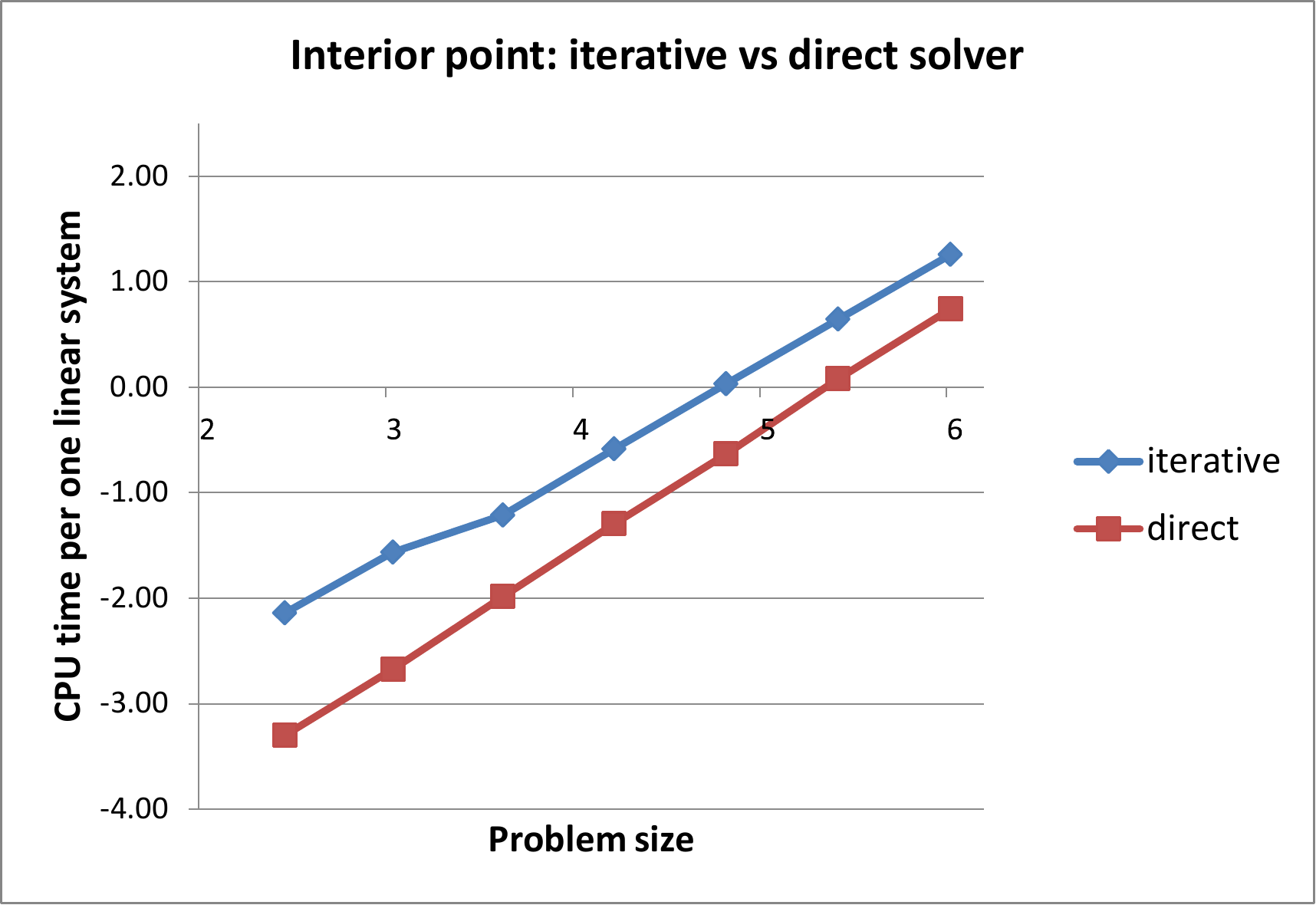}}
  \end{center}
  \caption{Example~2, \emph{left}: total CPU time spent in the iterative linear solver for the interior point and the OC method; \emph{right}: interior point method, total CPU time spent in the iterative linear solver and in the direct solver.}\label{fig:22}
\end{figure}

Figure~\ref{fig:22} (left) gives the comparison of the interior point with the
OC method. We can see even more clearly than in the previous example the faster
growth of the OC method. When we calculate the degree of the assumed polynomial
function $cn^d$ of the problem dimension $n$ from the larger examples, we will
obtain $d=0.944$ for the interior point method (so a linear growth) and
$d=1.28$ for the OC method.

%UPDATE UPDATE FIGURES

%%%%%%%%%%%%%%%%%%%%%%%%%%%%%%%%%%%%%%%%
%%%%%%%%%%%%%%%%%%%%%%%%%%%%%%%%%%%%%%%%

Figure~\ref{fig:22} (right) compares the iterative solver used in the interior
point method with the Cholesky solver (see the beginning of this section),
giving the same picture as in the previous example.

\begin{figure}[h]
\begin{center}
 \resizebox{0.48\hsize}{!}
   {\includegraphics{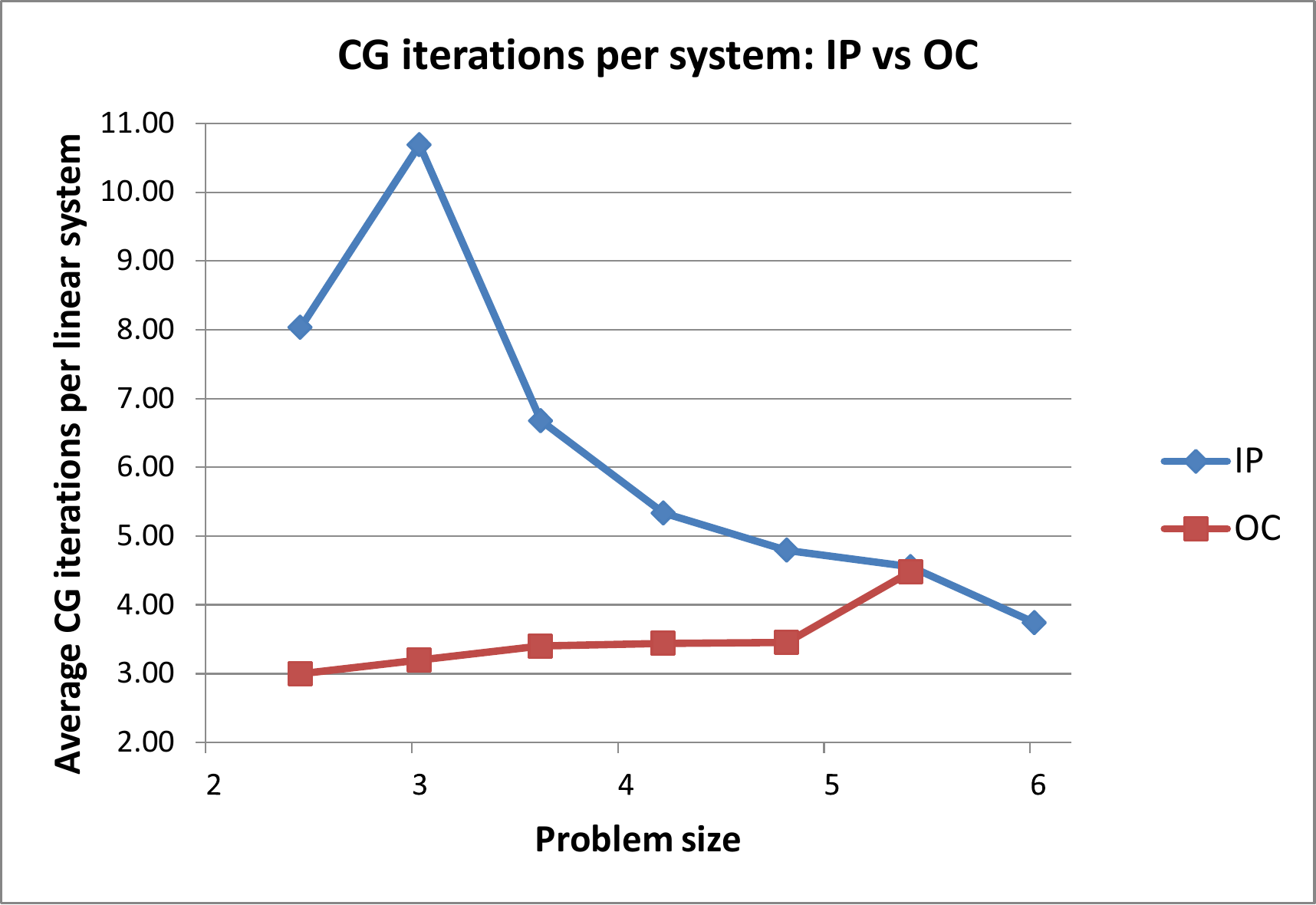}}
  \end{center}
  \caption{Example~2, average number of CG iterations per linear system for the interior point and the OC method.}\label{fig:24}
\end{figure}

Finally in Figure~\ref{fig:24} we compare the average number of CG steps per
linear system in the interior point and the OC solver. We can see that while
the graph is decreasing for the IP method, it is slowly increasing in case of
the OC method. The reason for that is that, in this example, we had to decrease
the stopping criterium for the CG solver in the OC method, in order to
guarantee its convergence (see Section~\ref{sec:CGOC} for explanation).
%%%%%%%%%%%%%%%%%%%%%%%%%%%%%%%%%%%%%%%%
%%%%%%%%%%%%%%%%%%%%%%%%%%%%%%%%%%%%%%%%

\subsection{Example 3}
\begin{figure}[h]
\begin{center}
 \resizebox{0.85\hsize}{!}
   {\includegraphics{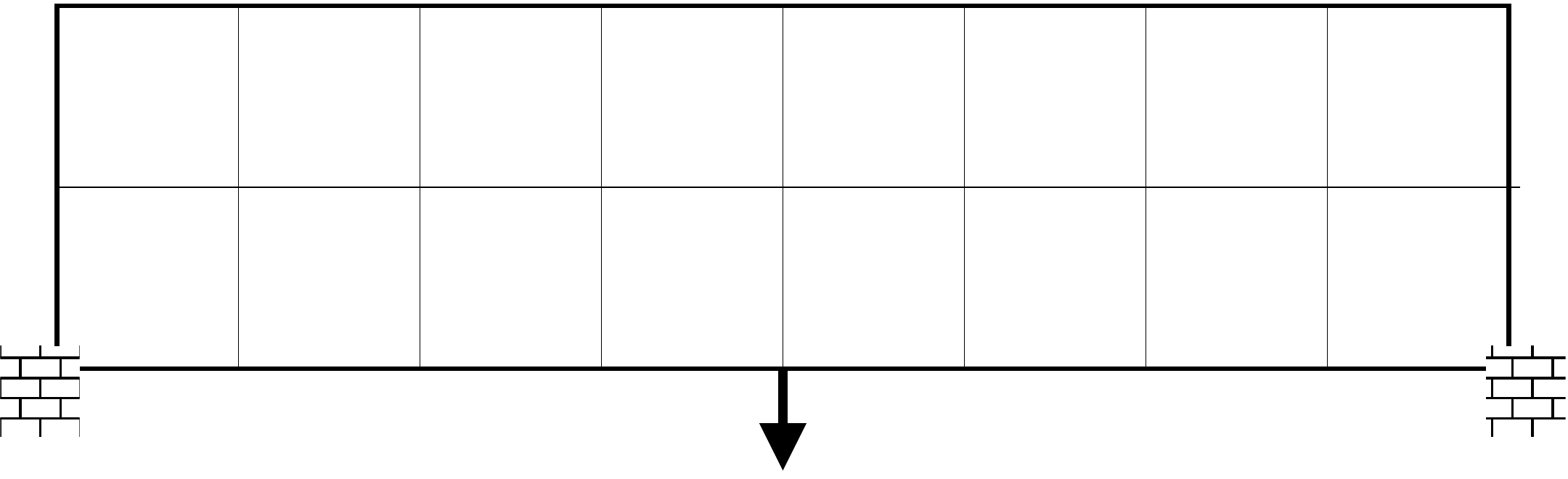}}\\
    \resizebox{0.8\hsize}{!}
   {\includegraphics{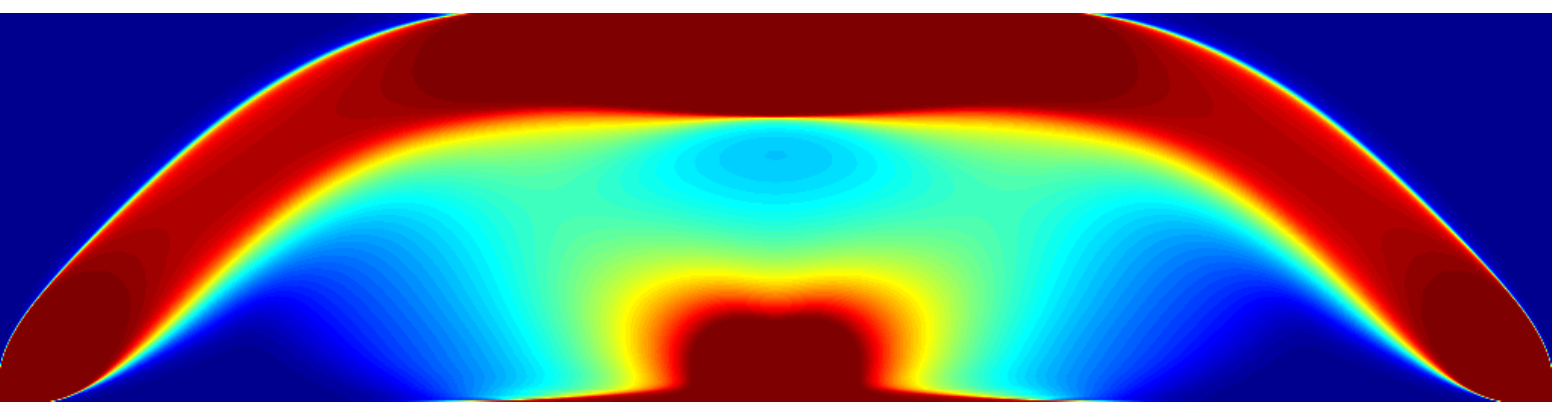}}
  \end{center}
  \caption{Example~3, initial setting with coarsest mesh and optimal solution.}\label{fig:31}
\end{figure}

The computational domain for our final example is a rectangle, initially
discretized by $8\times 2$ finite elements. The two corner points on the lower
edge are fixed and a vertical load is applied in the middle point of this edge;
see Figure~\ref{fig:31}. We use up to eight refinement levels with the finest
mesh having 262\,144 elements and 568\,850 nodal variables (after elimination
of the fixed nodes).
%%%%%%%%%%%%%%%%%%%%%%%%%%%%%%%%%%%%%%%%
%%%%%%%%%%%%%%%%%%%%%%%%%%%%%%%%%%%%%%%%

The results of the interior point method are shown in Table~\ref{tab:5}. Yet
again, the total number of CG iterations is decreasing with the increasing size
of the problem and the average number of CG iterations per linear system is
very low and not increasing. The negative complexity factor is caused by the
exceptional difficulties of the CG method in the last interior point step in
problem 823.
% Table generated by Excel2LaTeX from sheet 'br_school (2)'
\begin{table}[htbp]
  \centering
  \caption{Example 3, IP method with iterative solver}
    \begin{tabular}{crccccc}
    \toprule
              &       &       &  total     &   solver    & average  \\
    problem & variables  & feval & CG iters & CPU time & CG iters  \\
    \midrule
    822   & 170    & 33    & 284   & 0.21  & 8.61  \\
    823   & 594    & 31    & 383   & 0.78  & 12.35 \\
    824   & 2210   & 32    & 121   & 0.60  & 3.78  \\
    825   & 8514   & 31    & 166   & 3.41  & 5.35  \\
    826   & 33410  & 26    & 140   & 14.8  & 5.38  \\
    827   & 132354 & 26    & 133   & 78.8  & 5.12  \\
    828   & 526850 & 25    & 121   & 217.0 & 4.84  \\
    \bottomrule
    \end{tabular}%
  \label{tab:5}%
\end{table}%

Table~\ref{tab:3} presents the results of the OC method. As in Example~2, the
average number of CG iterations is increasing due to the worse conditioning.
% Table generated by Excel2LaTeX from sheet 'br_school (2)'
\begin{table}[htbp]
  \centering
  \caption{Example 3, OC method with iterative solver}
    \begin{tabular}{crccccc}
    \toprule
            &            &       &  total   &   solver & average   \\
    problem & variables  & feval & CG iters & CPU time & CG iters  \\
    \midrule
    822   & 170    & 23    & 69    & 0.04   & 3.00 \\
    823   & 594    & 37    & 147   & 0.21   & 3.97 \\
    824   & 2210   & 57    & 267   & 1.16   & 4.68 \\
    825   & 8514   & 75    & 374   & 7.40   & 4.99 \\
    826   & 33410  & 99    & 495   & 51.8   & 5.00 \\
    827   & 132354 & 111   & 665   & 290.0  & 5.99 \\
    828   & 526850 & 113   & 677   & 1250.0 & 5.99 \\
    \bottomrule
    \end{tabular}%
  \label{tab:6}%
\end{table}%

\begin{figure}[h]
\begin{center}
 \resizebox{0.5\hsize}{!}
   {\includegraphics{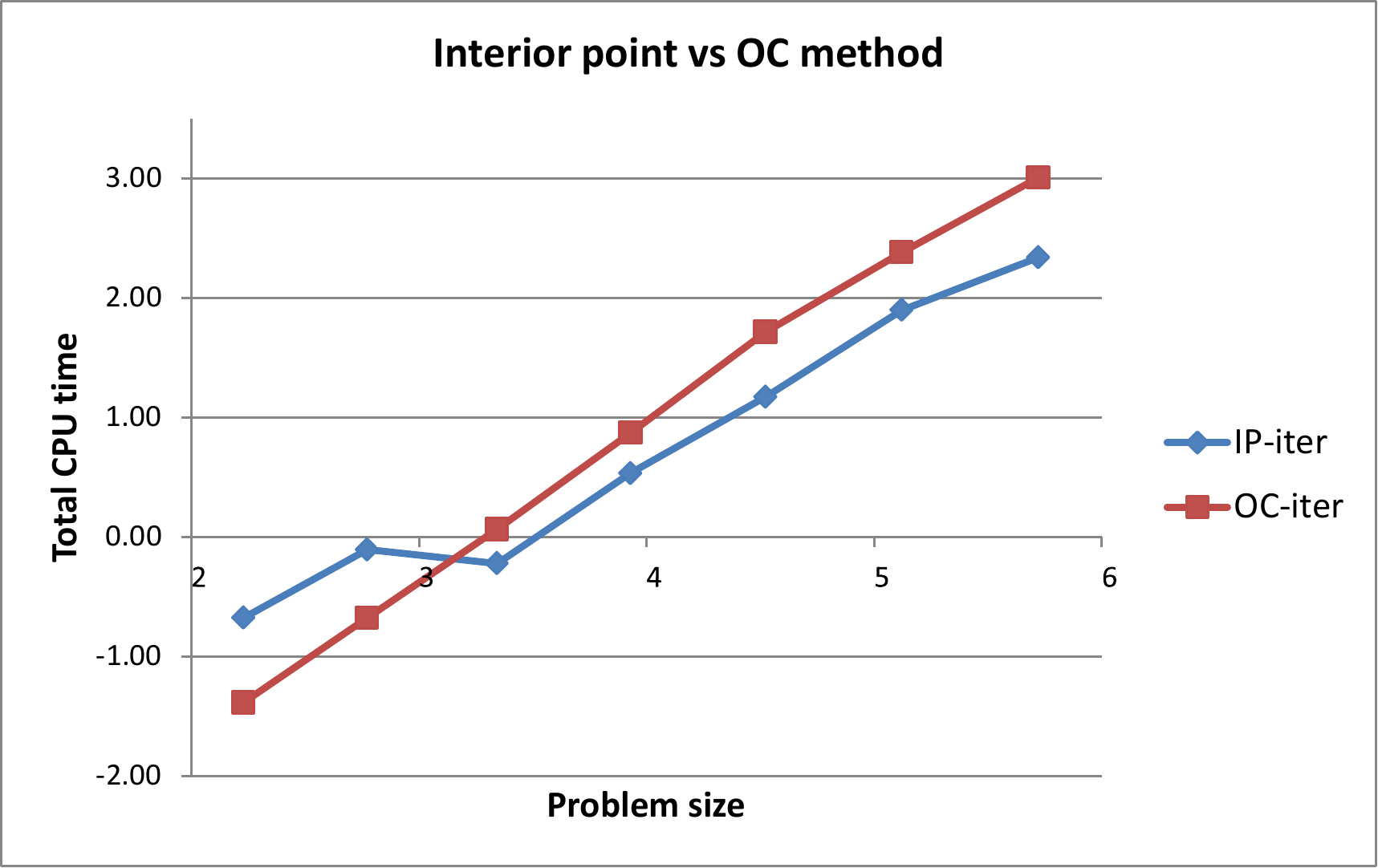}}\
 \resizebox{0.46\hsize}{!}
   {\includegraphics{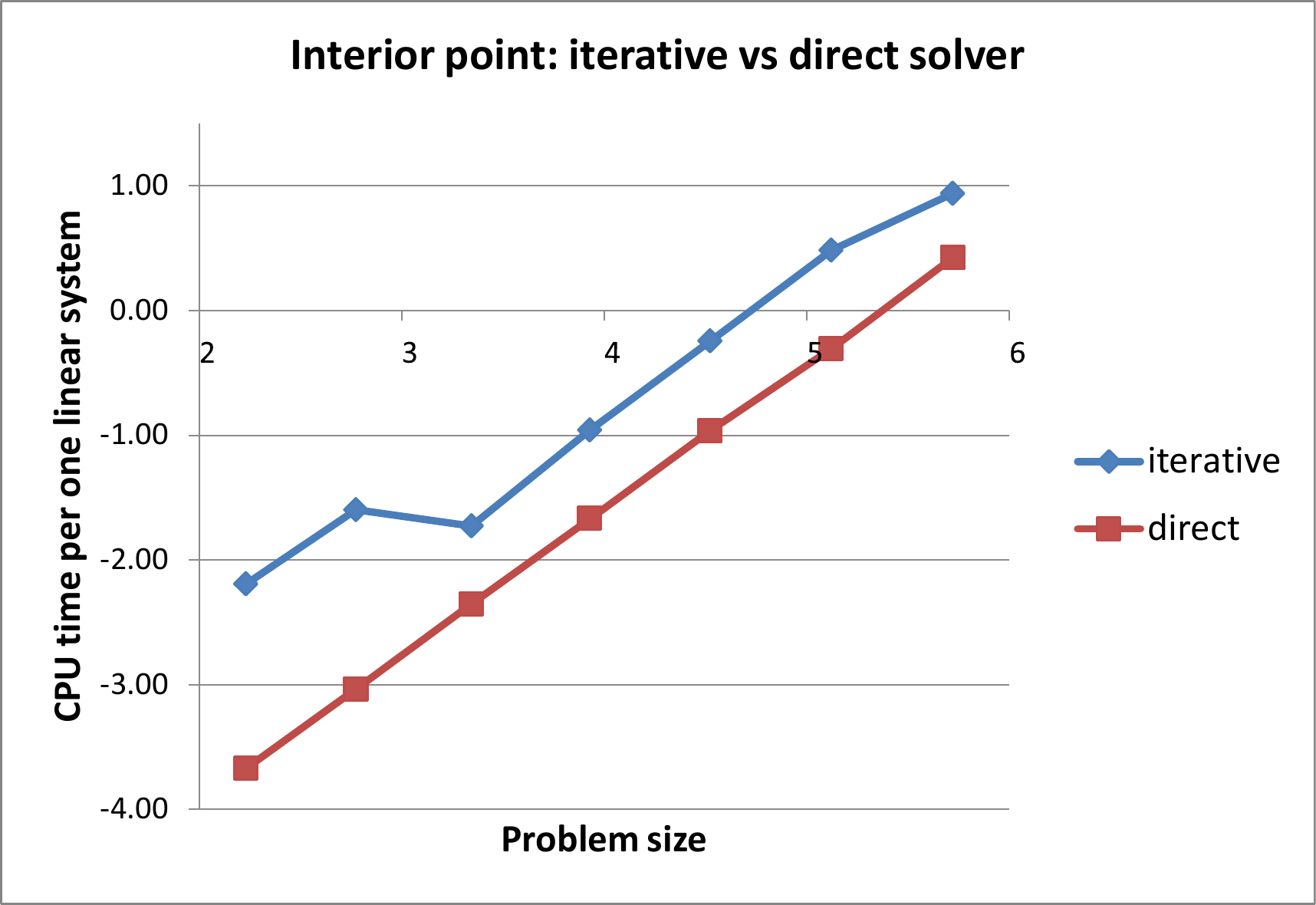}}
  \end{center}
  \caption{Example~3, \emph{left}: total CPU time spent in the iterative linear solver for the interior point and the OC method; \emph{right}: interior point method, total CPU time spent in the iterative linear solver and in the direct solver.}\label{fig:32}
\end{figure}

Figure~\ref{fig:32} (left) compares of the interior point with the OC method.
Yet again, the interior point method is a clear winner, both in the absolute
timing as in the growth tendency. Calculating the degree of the assumed
polynomial function $cn^d$ of the problem dimension $n$ from the larger
examples, we get $d=1.09$ for the interior point method and $d=1.24$ for the OC
method.
%%%%%%%%%%%%%%%%%%%%%%%%%%%%%%%%%%%%%%%%
%%%%%%%%%%%%%%%%%%%%%%%%%%%%%%%%%%%%%%%%

\begin{figure}[h]
\begin{center}
 \resizebox{0.48\hsize}{!}
   {\includegraphics{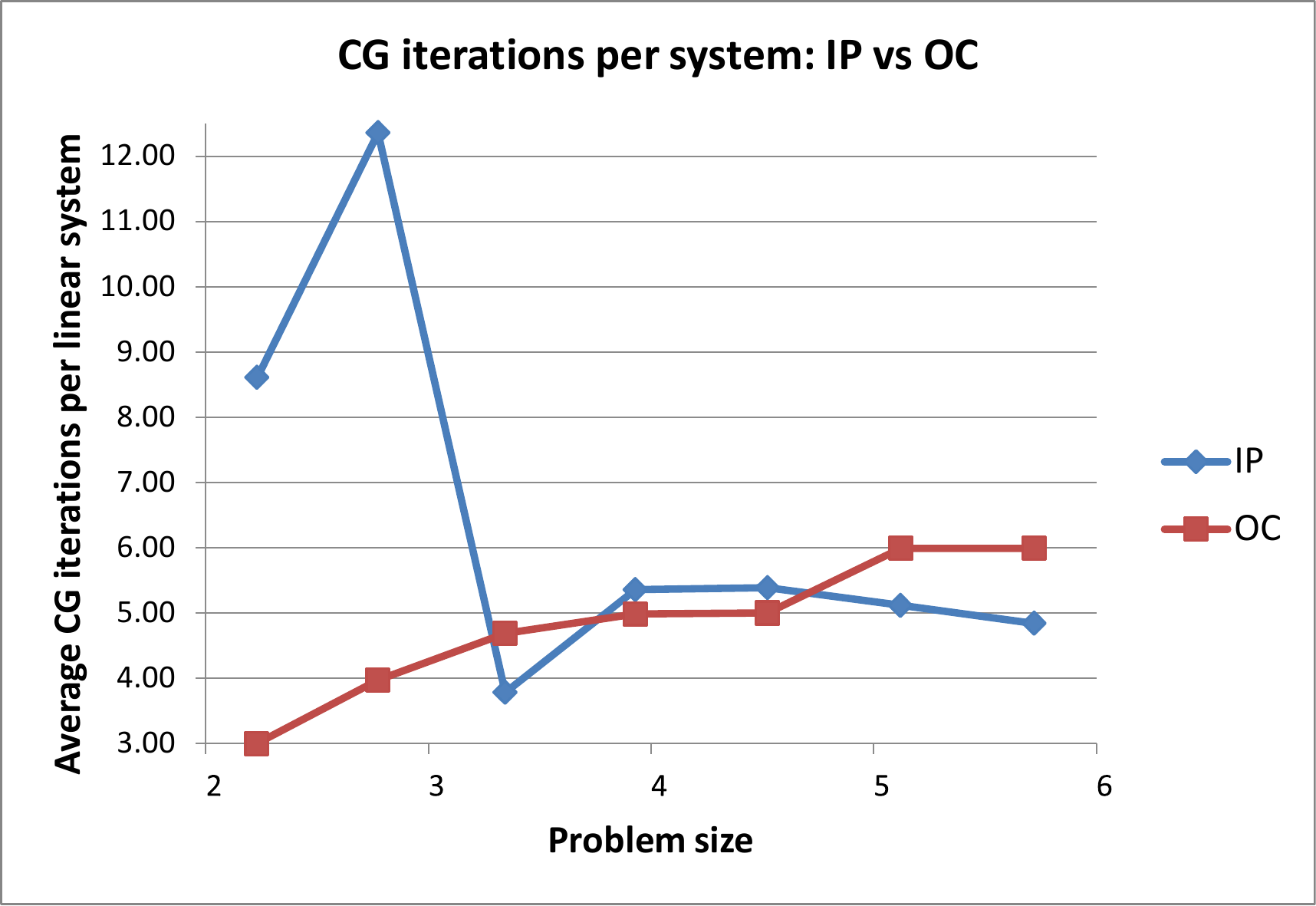}}
  \end{center}
  \caption{Example~3, average number of CG iterations per linear system for the interior point and the OC method.}\label{fig:34}
\end{figure}

In Figure~\ref{fig:32} (right) we compare the iterative solver used in the
interior point method with the Cholesky solver (see the beginning of this
section).
Finally in Figure~\ref{fig:34} we compare the average number of CG steps per
linear system in the interior point and the OC solver. We can see that while
the graph for the IP method has a decreasing tendency, it is increasing in case
of the OC method. As before, the reason for that is that we had to decrease the
stopping criterium for the CG solver, in order to guarantee its convergence
(see Section~\ref{sec:CGOC}).
%%%%%%%%%%%%%%%%%%%%%%%%%%%%%%%%%%%%%%%%
%%%%%%%%%%%%%%%%%%%%%%%%%%%%%%%%%%%%%%%%

\section{How exact is `exact'?}\label{sec:exact}
\subsection{Interior point method}\label{sec:exactIP}
\begin{figure}[h]
\begin{center}
 \resizebox{0.48\hsize}{!}
   {\includegraphics{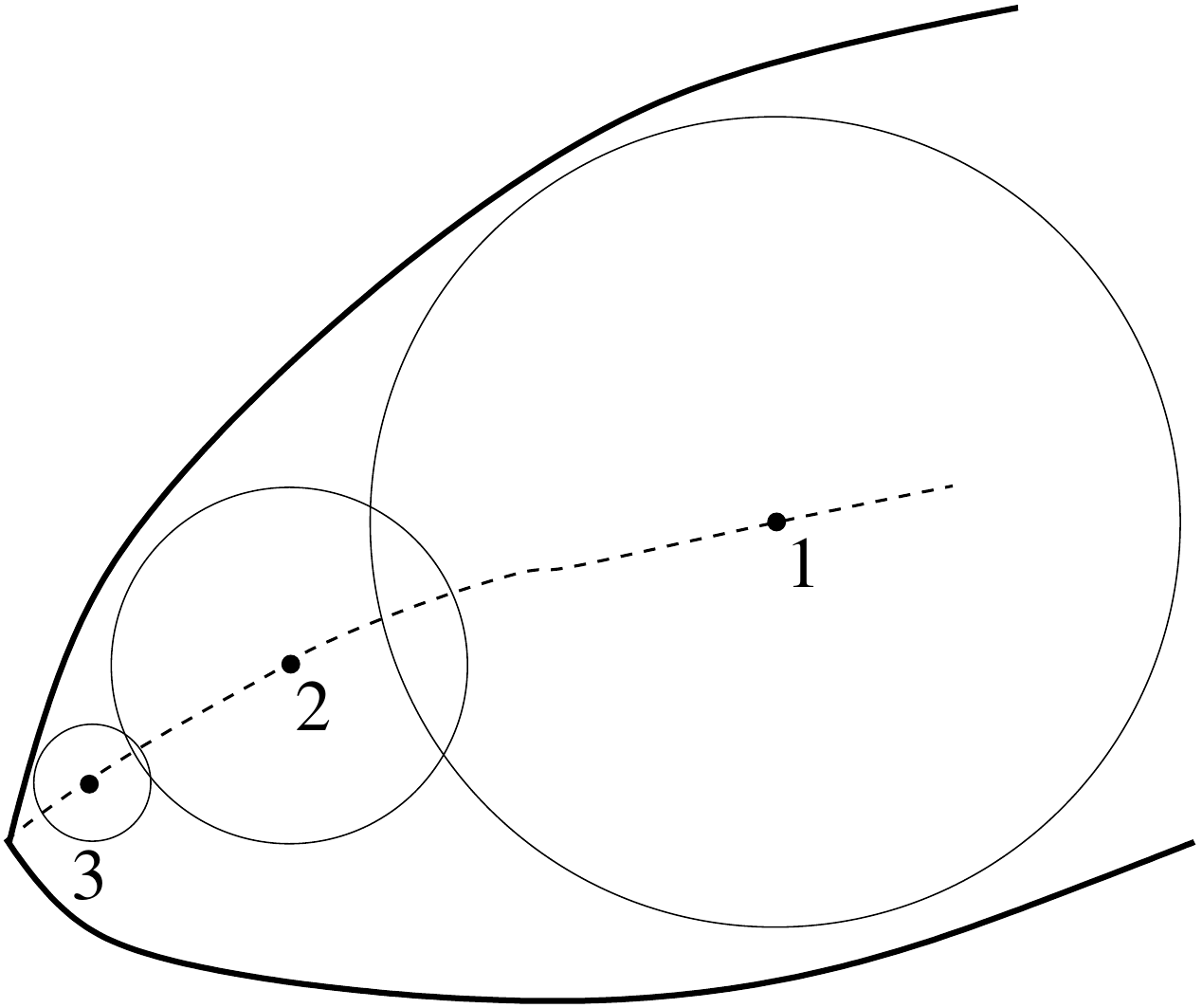}}
  \end{center}
  \caption{Illustration of interior-point method.}\label{fig:41}
\end{figure}

In this article we are using slightly nonstandard stopping criteria within the
interior point method. In particular, with the decreasing barrier parameters
$r,s$ we do \emph{not} decrease the stopping tolerances
$\tau_{\scriptscriptstyle\rm NWT}$ and $\tau_{\scriptscriptstyle\rm CG}$ for
the Newton method and for the conjugate gradients, respectively, although both
is required for the theoretical convergence proof. In Figure~\ref{fig:41} we
try to give a schematic explanation. Here we depict the feasible region and
three points $x_1,x_2,x_3$ on the central path, corresponding to three values
of the barrier parameter $r_1>r_2>r_3$. The exact solution lies in the corner
of the feasible region. The circle around each of these points depict the
region of stopping tolerance of the Newton method, once we get within, the
Newton method will stop. The radius of these circles is decreasing, even though
$\tau_{\scriptscriptstyle\rm NWT}$ is kept constant.
%%%%%%%%%%%%%%%%%%%%%%%%%%%%%%%%%%%%%%%%
%%%%%%%%%%%%%%%%%%%%%%%%%%%%%%%%%%%%%%%%
The idea is now obvious: it is ``better'' to stay within the tolerance circle
of $x_3$ rather than to get very close to $x_2$.

In the lemma below, $x^*$ is a point on the central path corresponding to a
barrier parameter $s$ and $x$ an approximation of $x^*$ resulting from inexact
Newton method. We will show that, even with a fixed stopping criterium for the
Newton method, $x$ must converge to $x^*$ with $s$ going to zero. For
simplicity of notation, we will just verify it for the lower bound
complementarity part of $\widetilde{\rm Res}^{(3)}$.
\begin{Lemma}
Let $x^*>0$ satisfies the perturbed scaled complementary condition
\begin{equation}\label{eq:111}
  \frac{\varphi_i x^*_i - s}{x_i} = 0,\quad i=1,\ldots,m
\end{equation}
and let $x>0$ be an approximation of $x^*$ satisfying
\begin{equation}\label{eq:112}
  \|z\|\leq \tau,\quad z_i=\frac{\varphi_i x_i - s}{x_i}
\end{equation}
with some $\tau>0$. Then there is an $\varepsilon>0$ depending on $s$ and
$\tau$ such that $\|x^*-x\|\leq\varepsilon$. Moreover, if $s$ tends to zero
then also $\varepsilon$ tends to zero.
\end{Lemma}
\begin{Proof}
From (\ref{eq:111}) we have that $\varphi_i = \frac{s}{x_i^*}$ and thus
(\ref{eq:112}) can be written as
$$
  \sum_{i=1}^m\left(\frac{s}{x_i^*} - \frac{s}{x_i}\right)^2\leq \tau^2
$$
which is, in particular, means that
$$
  \left|\frac{s}{x_i^*} - \frac{s}{x_i}\right| \leq \tau,\quad i=1,\ldots,m\,,
$$
i.e.,
$$
  \frac{|x_i^*-x_i|}{x^*_ix_i} \leq s\tau,\quad i=1,\ldots,m\,.
$$
Clearly, when $s$ tends to zero, $x$ must tend to $x^*$.
\end{Proof}

How good solution can we get when replacing the (``exact'') direct solver by an
inexact iterative method for the solution of the Newton systems? We may expect
that, with the ever decreasing barrier parameter, the inexact version will get
into numerical difficulties sooner than the exact one. Table~\ref{tab:IP}
answers this question. In topology optimization, the important variable is $x$,
the ``density''. With lower bound equal to zero, the quality of the solution
may be characterized by the closeness of components of $x$ to this lower bound
(that is, in examples where the lower bound is expected to be reached, such as
in Example~1 with sufficiently fine discretization). In Table~\ref{tab:IP} we
display the smallest component of $x$, denoted by $x_{\rm min}$ for Example~1
with 6 refinements levels, i.e., example 226 from Table~\ref{tab:2}. The
meaning of other columns in Table~\ref{tab:IP} is the following:
\begin{description}
\setlength\itemsep{0em}
\item[] {\bf barrier}\ldots{}the smallest value of the barrier parameters
    $s,r$ before the interior point algorithm was terminated;
\item[] {\bf IP,NWT,CG}\ldots{}the total number of iterations of the
    interior point method, the Newton method and conjugate gradients,
    respectively;
\item[] {\bf Cholesky}\ldots{}the linear system was solved by the CHOLMOD
    implementation of the Cholesky method;
\item[] {\bf CG tol fixed}\ldots{} the linear system was solved by the
    multigrid preconditioned conjugate gradient method with a fixed
    stopping criterium $\|r\|\,\|b\|\leq 10^{-2}$; see (\ref{eq:cgstop});
\item[] {\bf CG tol decreasing}\ldots{} as above but with a variable
    stopping criterium $\|r\|\,\|b\|\leq\tau_{\scriptscriptstyle\rm CG}$,
    where $\tau_{\scriptscriptstyle\rm CG}$ is initially equal to $10^{-2}$
    and is then multiplied by 0.5 after each major iteration of the
    interior point method.
\end{description}
\begin{table}[htbp]
  \centering
  \caption{Number of iterations and error in the IP solution for different values of $\tau_{\scriptscriptstyle\rm IP}$ and three different linear
  solvers.}
    \begin{tabular}{rr|rc|rrc|rrc}
    \toprule
     &       & \multicolumn{2}{c|}{Cholesky} & \multicolumn{3}{c|}{CG tol fixed} & \multicolumn{3}{c}{CG tol decreasing} \\
    \midrule
    $\tau_{\scriptscriptstyle\rm IP}$ & IP  & NWT  & $x_{\rm min}$ & NWT  & CG  & $x_{\rm min}$ & NWT  & CG  & $x_{\rm min}$ \\\midrule
    $10^{-8}$  & 12 & 28 & $1.6\cdot 10^{-5}$  & 28  & 139   & $1.8\cdot 10^{-5}$ & 28 & 587   & $1.6\cdot 10^{-5}$\\
    $10^{-10}$ & 15 & 34 & $1.0\cdot 10^{-7}$  & 35  & 291   & $2.7\cdot 10^{-7}$ & 34 & 4285  & $1.0\cdot 10^{-7}$\\
    $10^{-12}$ & 18 & 40 & $1.0\cdot 10^{-9}$  & 72  & 2832  & $2.4\cdot 10^{-9}$ & 40 & 10042 & $1.5\cdot 10^{-9}$\\
    $10^{-14}$ & 21 & 46 & $6.4\cdot 10^{-12}$ & 296 & 63674 & $1.9\cdot 10^{-11}$& 53 & 23042 & $1.9\cdot 10^{-11}$\\
    $10^{-16}$ & 24 & 52 & $6.2\cdot 10^{-14}$ & 489 & 88684 & $1.4\cdot 10^{-13}$& 82 & 52042 & $1.4\cdot 10^{-13}$\\
    \bottomrule
    \end{tabular}%
  \label{tab:IP}%
\end{table}%
We can see that all three algorithms were able to solve the problem to very
high accuracy. However, both versions of the CG method had problems with very
low values of the barrier parameter. The ``CG tol fixed'' version needed very
high number of the Newton steps, while the ``CG tol decreasing'' version needed
very high number of the CG steps to reach the increased accuracy. (Notice that
the maximum number of CG iterations for one system was limited to 1000.) On the
other hand, for barrier parameter equal to $10^{-8}$ (our choice in the
numerical examples above), both inexact solvers were on par with the exact one
and, due to the lower accuracy required and thus lower number of CG steps, the
``CG tol fixed'' version is the method of choice.

\subsection{OC method}
In the OC method, we have to solve the equilibrium problem with the stiffness
matrix $K(x)$; that means, $K(x)$ must not be singular. A common way how to
approach this is to assume that $x$ is strictly positive, though very small.
Typically, one would modify the lower bound constraint to $0<\underline{x}\leq
x_i$, $i=1,\ldots,m$ with $\underline{x}=10^{-6}$, for instance. Once the OC
method is terminated, all values of $x$ with $x_i=\underline{x}$ are set to
zero. This is usually considered a weakness of the OC method, because we do not
exactly solve the original problem, only its approximation (see
\cite{achtziger}). Somewhat surprisingly, in the examples we solved using our
MATLAB code, the value of $\underline{x}$ could be actually very low, such as
$\underline{x}=10^{-30}$. The stiffness matrix $K(x)$ will, consequently,
become extremely ill-conditioned (in the above case the condition number will
be of order $10^{30}$), nevertheless, CHOLMOD does not seem to have a problem
with that and the OC method converges in about the same number of iterations as
if we set $\underline{x}=10^{-6}$.

The main question is how does the quality of the solution depends on the
heuristic stopping criterium (\ref{eq:OCstop}). Our next Table~\ref{tab:exact1}
sheds some light on this. We solve the example 226 from Table~\ref{tab:2} for
various values of the stopping criterium $\tau_{\scriptscriptstyle\rm OC}$ and
two different values of the lower bound $\underline{x}$. We then compute,
pair-wise, the norm of the difference of these solution. The notation
$\tau_{\scriptscriptstyle\rm OC}=10^{-\inf}$ is used for the case when the
stopping criterion (\ref{eq:OCstop}) is ignored and the OC method is terminated
after a very high number of iterations, in this case 5000 (i.e., 10000
solutions of the linear system). The resulting solution serves as the best
approximation of the exact solution that can be obtained within the
computational framework. So looking at Table~\ref{tab:exact1}, we can see that,
for instance, the maximum norm of the difference between the solutions with
$\tau_{\scriptscriptstyle\rm OC}=10^{-5}$ and $\tau_{\scriptscriptstyle\rm
OC}=10^{-\inf}$ is $\|x_{-5}-x_{-\inf}\|_\infty=0.126$, while
$\|x_{-9}-x_{-\inf}\|_\infty=0.016$. Notice that the norm is not scaled, e.g.,
by the dimension of $x$, hence the numbers are relatively large. Also, to get a
clearer picture, we used a direct linear system solver.
% Table generated by Excel2LaTeX from sheet 'Sheet3'
\begin{table}[htbp]
  \centering
  \caption{The norm of difference of two OC solutions $x$ for various values of the stopping criterium $\tau_{\scriptscriptstyle\rm OC}=10^{-5},10^{-7},10^{-9},10^{-\inf}$,
  and for two values of the lower bound $\underline{x}=10^{-7}$ and $\underline{x}=10^{-17}$. Upper triangle shows the 2-norm, lower triangle the infinity norm.}
    \begin{tabular}{r|rrrrr|rrrr}
    \toprule
           lower bound&       & \multicolumn{4}{c|}{$10^{-7}$}  & \multicolumn{4}{c}{$10^{-17}$} \\ \midrule
                 & \multicolumn{1}{|c|}{$\tau_{\scriptscriptstyle\rm OC}$} & -5 & -7 & -9 & ${-\inf}$ & -5 & -7 & -9 & ${-\inf}$ \\ \midrule
     \multicolumn{1}{c|}{\multirow{4}[2]{*}{$10^{-7}$}} & \multicolumn{1}{c|}{-5} & 0     & 1.09  & 1.19  & 1.26  & $2\cdot 10^{-6}$ & 1.10   & 1.18  & 1.26 \\
     \multicolumn{1}{c|}{} & \multicolumn{1}{c|}{-7}& 0.114 & 0     & 0.116 & 0.281 & 1.09  & 0.013 & 0.110  & 0.281 \\
     \multicolumn{1}{c|}{} & \multicolumn{1}{c|}{-9} & 0.123 & 0.01 & 0     & 0.190 & 1.19  & 0.103 & 0.009 & 0.190 \\
    \multicolumn{1}{c|}{} & \multicolumn{1}{c|}{${-\inf}$} & 0.126 & 0.025 & 0.016 & 0     & 1.26  & 0.271 & 0.198 & $4\cdot 10^{-6}$ \\\midrule
    \multicolumn{1}{c|}{\multirow{4}[2]{*}{$10^{-17}$}} & \multicolumn{1}{c|}{-5} & 2e-7 & 0.114 & 0.123 & 0.126 & 0     & 1.10   & 1.19  & 1.26 \\
     \multicolumn{1}{c|}{} & \multicolumn{1}{c|}{-7} & 0.116 & 0.001 & 0.009 & 0.024 & 0.116 & 0     & 0.094 & 0.271 \\
    \multicolumn{1}{c|}{} & \multicolumn{1}{c|}{-9} & 0.123 & 0.009 & $8\cdot 10^{-4}$ & 0.017 & 0.123 & 0.008 & 0     & 0.198 \\
     \multicolumn{1}{c|}{} & \multicolumn{1}{c|}{${-\inf}$} & 0.126 & 0.025 & 0.016 & $3\cdot 10^{-7}$ & 0.126 & 0.024 & 0.017 & 0 \\
    \bottomrule
    \end{tabular}%
  \label{tab:exact1}%
\end{table}%

\subsection{Interior point versus OC method}
We again solve example 226 from Table~\ref{tab:2}, this time by the interior
point method with an exact linear solver and various stopping parameters
$\tau_{\scriptscriptstyle\rm IP}$. In Table~\ref{tab:exact2}, these solutions
are compared (in two different norms), to the `exact' solution obtained in the
previous section by 5000 iterations of the OC method with
$\underline{x}=10^{-17}$. Comparing these numbers to those in
Table~\ref{tab:exact1}, we can see that the IP method delivers very good
solution already for our standard value $\tau_{\scriptscriptstyle\rm
IP}=10^{-8}$; this is comparable to OC solution with
$\tau_{\scriptscriptstyle\rm OC}=10^{-7}$. Moreover, decrease of
$\tau_{\scriptscriptstyle\rm IP}$ leads to a rapid decrease of the error,
unlike in the OC method.
\begin{table}[htbp]
  \centering
  \caption{Two different norms of the error of the IP method in variable $x$ for different values of the stopping parameter $\tau_{\scriptscriptstyle\rm IP}$. As an `exact' solution $x^*$ we take the OC solution after 5000 iterations with lower bound $\underline{x}=10^{-17}$.}
    \begin{tabular}{ccc}
    \toprule
             $\tau_{\scriptscriptstyle\rm IP}$    & $\|x-x^*\|_2$ & $\|x-x^*\|_{\infty}$ \\\midrule
     $10^{-8}$ & $2.47\cdot 10^{-1}$ & $2.49\cdot 10^{-2}$ \\
     $10^{-10}$ & $9.60\cdot 10^{-3}$ & $1.50\cdot 10^{-3}$ \\
     $10^{-12}$ & $2.07\cdot 10^{-4}$ & $4.95\cdot 10^{-5}$ \\
     $10^{-14}$ & $1.70\cdot 10^{-6}$ & $4.66\cdot 10^{-7}$ \\
    \bottomrule
    \end{tabular}%
  \label{tab:exact2}%
\end{table}%

\section{The SIMP model of topology optimization}
\begin{figure}[h]
\begin{center}
 \resizebox{0.48\hsize}{!}
   {\includegraphics{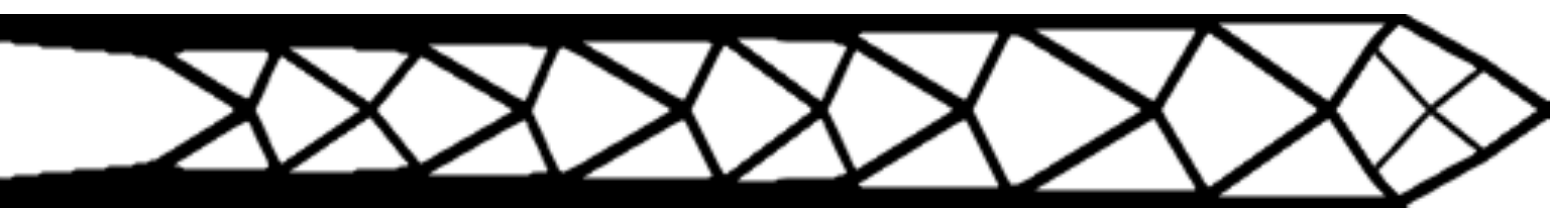}}
 \resizebox{0.48\hsize}{!}
   {\includegraphics{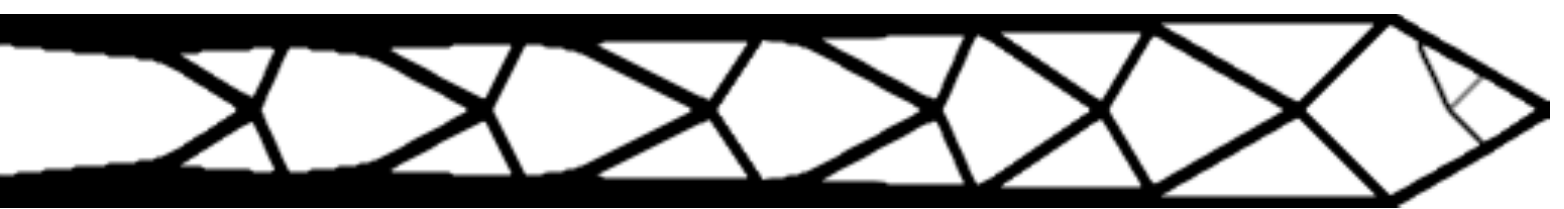}}\\
 \resizebox{0.48\hsize}{!}
   {\includegraphics{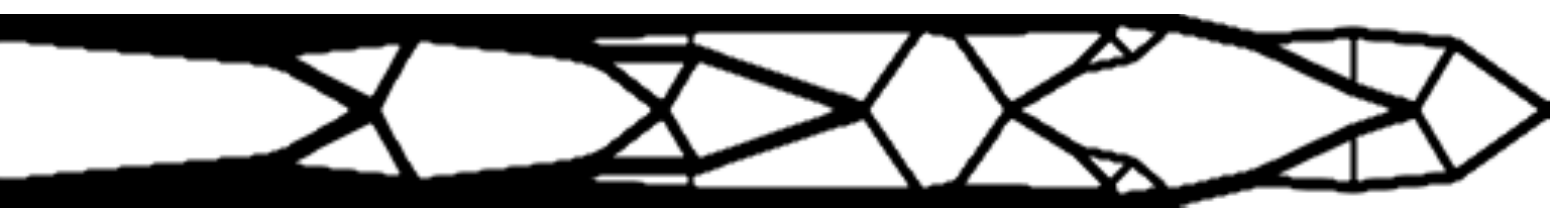}}
 \resizebox{0.48\hsize}{!}
   {\includegraphics{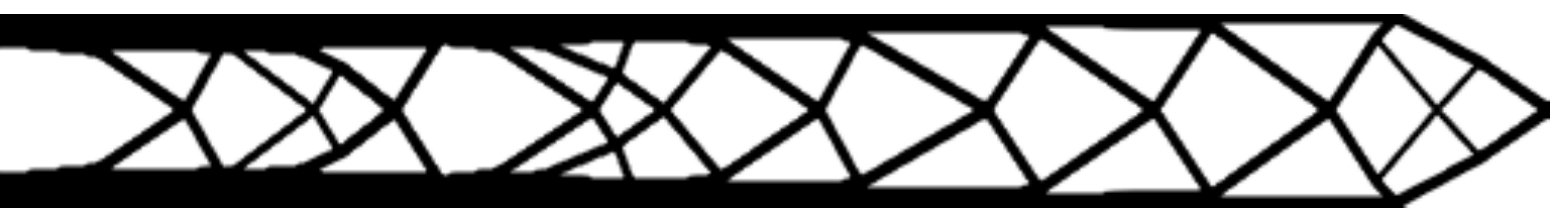}}
  \end{center}
  \caption{Optimal result SIMP of code {\tt top88} starting from the default
  initial point and from three other initial points.}\label{fig:101}
\end{figure}

A natural question arises about the applicability of the presented approach to
a more popular model of topology optimization, namely, the Solid Isotropic
Material with Penalisation (SIMP) \cite{bendsoe-sigmund}. When used with a
suitable filtering, one can guarantee at least the existence of a solution of
the infinite-dimensional problem and convergence of the finite element
discretization to this solution \cite{bourdin2001filters}. However, the problem
is non-convex and exhibits many local minima, as demonstrated in
Figure~\ref{fig:101}. There we show solutions obtained by the code {\tt top88}
\cite{top88} from various initial points (with a load vector modified to the
original code). The calling sequence of the code was {\tt
top88(512,64,1,3,1.5,2)}. The ordering of the plots in Figure~\ref{fig:101} is
$\small\begin{bmatrix}(a) &(b)\\(c)&(d)
\end{bmatrix}$ and the respective values of the computed optimal compliance
are (a) 128.402 (from the default starting point), (b) 126.773, (c) 134.359,
(d) 129.002. Notice that we strengthened the stopping criterion of {\tt top88}
from $10^{-2}$ to $10^{-3}$.
%
%%%%%%%%%%%%%%%%%%%%%%%%%%%%%%%%%%%%%%%%
%%%%%%%%%%%%%%%%%%%%%%%%%%%%%%%%%%%%%%%%

It is thus difficult to compare the efficiency of various algorithms, as each
may converge to a different local solution (see also \cite{rojas-stolpe}).
Moreover, starting from two different initial points, an algorithm may converge
to two different solutions and the number of iterations needed to find the
respective solutions can differ substantially. For these reasons, we only give
a brief comparison of the IP method with an exact solver and with an iterative
solver, demonstrating that the iterative method is still a viable and efficient
option. (A similar comparison for the OC method can be found in \cite{amir}.)

The SIMP model with the so-called density filter consists in a modification of
our original problem (\ref{eq:to}) where we replace the equilibrium equation $
  K(x) u = f
$ by
$$
  \widehat{K}(\tilde{x}) u = f\qquad\mbox{with} \qquad \widehat{K}(\tilde{x}) = \sum_{i=1}^m \tilde{x}_i^p K_i\,,
$$
where $\tilde{x}_i$ is computed as a weighted average of the values of $x$ in a
close neighborhood of the $i$-th element. More precisely,
$$
  \tilde{x} = W x \quad{\rm with}\quad
  W_{:i} = \frac{1}{\sum_{j=1}^m \widehat{W}_{ij}} \widehat{W}_{:i}
\quad{\rm and}\quad
  \widehat{W}_{ij} = \max(0,r_{\rm min} - {\rm dist}(i,j))\,,
$$
where $W_{ij}$ is the $(i,j)$-th element of matrix $W$ and $W_{:i}$ denotes the
$i$-th column of $W$. Here ${\rm dist}(i,j)$ is a function measuring the
distance between the $i$-th and $j$-th element (e.g.\ Manhattan distance or
Euclidean distance of element centers), and $r_{\rm min}$ is a given radius of
the density filter. The typical choice of $p$ is $p=3$.

The interior point method from Sections 2 and 3 has to be adapted to the SIMP
model. In particular, the KKT condition (\ref{eq:KKT3}) will change to
$$
  -\frac{1}{2}u^T \left[ p\sum_{j=1}^m W_{ji}(Wx)_j^{p-1}K_j\right] u
  - \lambda - \varphi_i + \psi_i = 0,\quad i=1,\ldots,m
$$
which means that the matrix $B(u)$ in the linear system (\ref{eq:nwtr}) (and
thus (\ref{eq:nwtr2})) will be replaced by
$$
  \widehat{B}(u)=\left(p\sum_{j=1}^m W_{j1}(Wx)_j^{p-1}K_j u, \ldots,
  p\sum_{j=1}^m W_{jm}(Wx)_j^{p-1}K_j u\right)\,.
$$
The consequence of using the density filter is that the matrix
$\widehat{B}(u)\widehat{B}(u)^T$ (and thus the matrix $Z$ in (\ref{eq:nwtr2}))
will have bigger band-width and thus the Cholesky factors will have more
non-zero elements.

The purpose of the next example is merely to show that the iterative solver is
still a viable alternative to the direct one, even for the non-convex problem.
The reader should not be too concerned with the behaviour of the interior point
method, as it still uses the vanilla algorithm and explicit choice of step
length that were suitable for the convex problem but should be upgraded for the
non-convex one. However, this was not the goal of this paper.

\subsection{Example 4}
\begin{figure}[h]
\begin{center}
 \resizebox{0.48\hsize}{!}
   {\includegraphics{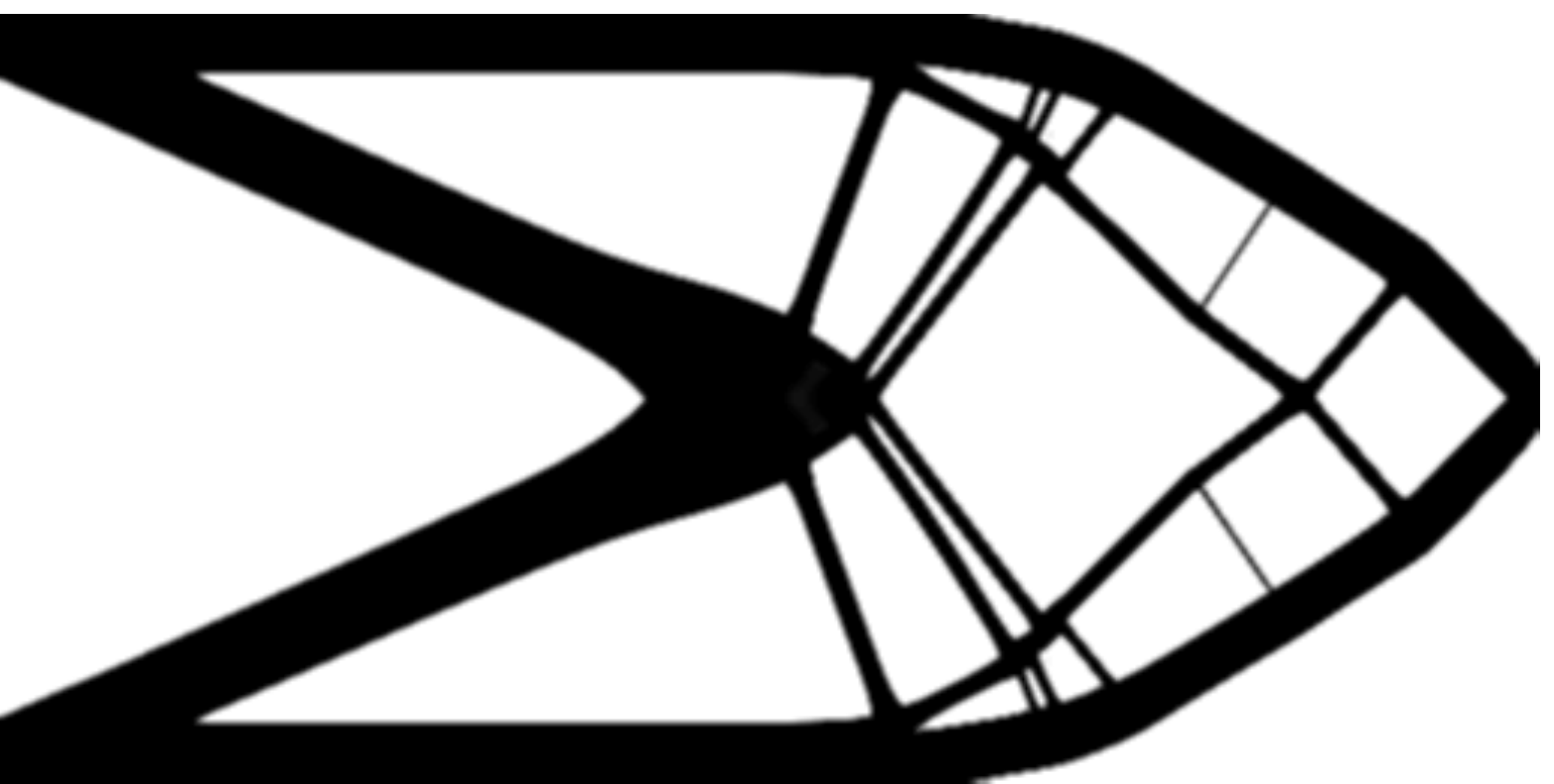}}
 \resizebox{0.48\hsize}{!}
   {\includegraphics{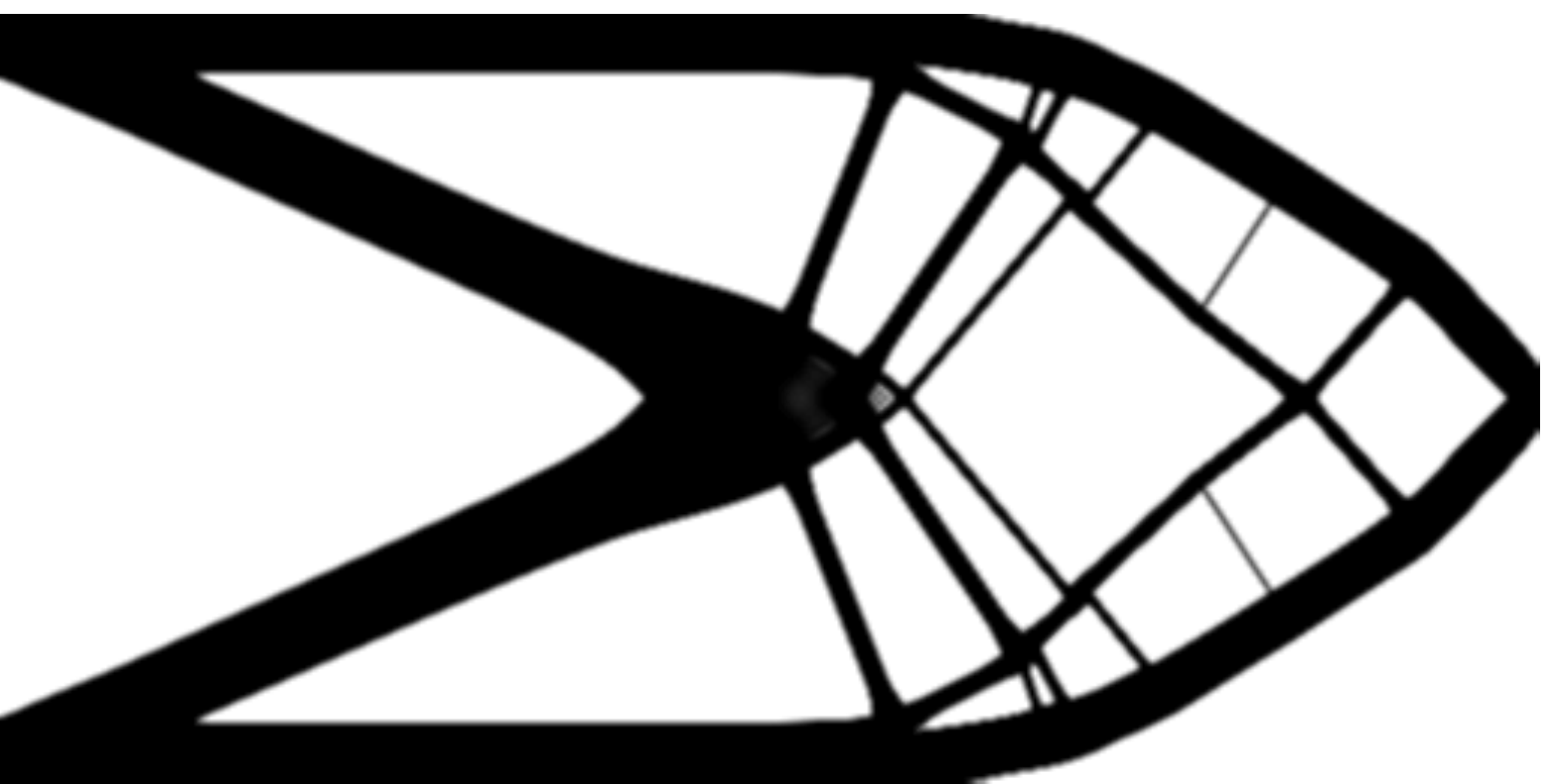}}
  \end{center}
  \caption{Example~4, optimal solution $\tilde{x}^*$ using the iterative (left)
  and the Cholesky (right) solver.}\label{fig:102}
\end{figure}

Consider again the problem from Example~2 with an upper bound $\overline{x}=3$,
this time solved using the SIMP model with density filter. The filter uses
Manhattan distance with $r_{\rm min}=2$. Figure~\ref{fig:102} shows the optimal
results for 8 refinement levels (problem {\tt 428}) when using the iterative
solver (left) and the Cholesky method (right). We can see that also in this
example the two versions of the code converge to different local minima with
almost identical objective function values (see Table~\ref{tab:201}).

%%%%%%%%%%%%%%%%%%%%%%%%%%%%%%%%%%%%%%%%
%%%%%%%%%%%%%%%%%%%%%%%%%%%%%%%%%%%%%%%%

The numbers presented in Table~\ref{tab:201} for problems {\tt 427} and {\tt
428} show the total number of Newton steps (feval), the CPU time only needed by
the linear solver and the optimal objective function value (obj). The number of
Newton steps needed by the IP method is not significantly influenced by the
choice of the solver. Also the number of CG iterations is still kept very low
with average 3.5 per linear system. To demonstrate that was the purpose of this
example.
\begin{table}[htbp]
  \centering
  \caption{Example 4, IP method with iterative and Cholesky solver}
    \begin{tabular}{cr|cccc|ccc}
    \toprule
              & &  \multicolumn{4}{c|}{CG} & \multicolumn{3}{c}{Cholesky}\\
    problem & variables  & feval & CG iters & time & obj & feval & time & obj  \\
    \midrule
    427   & 66\,048     & 151    & 514   &  316 & 3.5659 &  128 &  146 & 3.5681 \\
    428   & 263\,168    & 262    & 934   & 2410 & 3.4709 &  231 & 1500 & 3.4762\\
    \bottomrule
    \end{tabular}%
  \label{tab:201}%
\end{table}%

Finally, in Table~\ref{tab:202} we present the number of Newton steps needed in
every major IP iteration. The first row shows the value of the barrier
parameter $s$ (reduced in every iteration by factor 0.2). The next two rows
refer to problem {\tt 427} and show the number of Newton steps first when using
Cholesky method and then for the iterative solver. The final two rows are for
problem {\tt 428}. As we can see, most effort is spent in the last iterations,
unlike in the convex case, where the number of Newton steps was almost
constant. As mentioned before, a more sophisticated version of the IP method
would be needed for the non-convex case to avoid this behaviour.
\begin{table}[htbp]
  \centering
  \caption{Example 4, the barrier parameter $s$ and the corresponding number of
  Newton steps needed at every major IP iteration
  with the Cholesky and the CG solver, respectively.}
    \begin{tabular}{lccccccccccccc}
    \toprule
 &$s$& $5^0$ & $5^{-1}$ &  $5^{-2}$ & $5^{-3}$ & $5^{-4}$ & $5^{-5}$ & $5^{-6}$ & $5^{-7}$ & $5^{-8}$ & $5^{-9}$ & $5^{-10}$ & $5^{-11}$ \\
 \midrule
 \multirow{2}{1em}{427}&Chol&2&2&2&4&4&6&9&11&24&22&21&21\\
 &CG     &2&2&3&3&4&6&9&12&24&23&25&38\\
 \midrule
 \multirow{2}{1em}{428}&Chol&2&2&2&3&3&4&7&11&21&44&65&67\\
 &CG  &2&2&2&3&3&4&7&11&21&53&68&86\\
    \bottomrule
    \end{tabular}%
  \label{tab:202}%
\end{table}%

\section{Conclusions}
Based on the results of our numerical experiments, we make the following
conclusions. These conclusions only concern the convex problem.
\begin{itemize}
\item The interior point method clearly outperforms the OC method on
    large-scale problems. The larger the problem, the bigger the
    difference. This is independent of the fact whether direct or iterative
    solver is used for the linear system. It is also independent of the
    fact whether the linear systems (in both methods) are solve exactly or
    inexactly.
\item The inexact multigrid preconditioned CG method outperforms even a
    very sophisticated direct solver, at least for large to very-large
    scale problems. This holds for both, the interior point and the OC
    method.
\item The behaviour of the interior point method is very predictable. More
    surprisingly, also the behaviour of the chosen iterative method, the
    multigrid preconditioned conjugate gradients, is also very predictable
    and \emph{independent on the size of the problem}.
\item Also in the OC method, the multigrid preconditioned CG algorithm is
    predictable and very stable, both with respect to the size of the
    problem and of the OC iteration (and thus of the condition number of
    the stiffness matrix). Perhaps rather surprisingly, not more than 10 CG
    iterations are needed, even when high precision of the OC method is
    required. This is the effect of the multigrid preconditioner: notice
    that in \cite{wang} the authors report about 100--200 CG steps needed
    (with a different preconditioner) and thus propose to use so-called
    recycling of the Krylov subspaces, in order to accelerate CG
    convergence speed. This is just not needed here, given the very low
    number of CG steps.
\item The OC method has one noticeable advantage over the interior point
    method. It can quickly identify the ``very strongly'' active
    constraints, those with large Lagrangian multiplier. Due to the
    projection of variables on the feasible set, the active variables are
    then exactly equal to the bounds. Contrary to that, the interior point
    method only approaches the boundary. This may be particularly
    significant in case of lower bounds, when the user has to decide which
    values are cut off and considered zero (and thus interpreted as void).
    Clearly, the lower bound for the OC method has to be positive but it
    can be set very low (e.g., $10^{-17}$) and is then exactly reached.
\end{itemize}

From the above, it seems to be obvious to recommend the interior point method
with multigrid preconditioned CG solver as the method of choice for large scale
topology optimization problem. However, we should keep in mind that the use of
multigrid is rather restricted by the assumed existence of regularly refined
finite element meshes. This is easily accomplished when using ``academic''
examples with regular computational domains such as squares, rectangles, prisms
and unions of these. For geometrically complex domains appearing in practical
examples, multigrid may not be so suitable or may even be unusable. In these
cases, we can resort to domain decomposition preconditioners. In
\cite{turner-kocvara-loghin} it was shown that, in connection with the interior
point method, they also lead to very efficient techniques for topology
optimization problems.

\bibliographystyle{plain}
\bibliography{topo_mgm_paper}
%----------------------------------------------------------------
\end{document}